\documentclass[11pt,a4paper]{article}

\usepackage{algorithm}
\usepackage{amsfonts}
\usepackage{amsmath}
\usepackage{amssymb}
\usepackage[noend]{algpseudocode}
\usepackage[titletoc]{appendix}
\usepackage{bbm}
\usepackage{booktabs}
\usepackage{color}
\usepackage[most]{tcolorbox}
\usepackage[inline]{enumitem}
\usepackage{xcolor}
\usepackage{eqparbox}
\usepackage{fancyhdr,lipsum}
\usepackage[symbol]{footmisc}
\usepackage{graphicx}
\usepackage{lineno}
\usepackage{longtable}
\usepackage{lscape}
\usepackage{makecell}
\usepackage{mathrsfs}
\usepackage{mathtools}
\usepackage{multirow}
\usepackage[numbers]{natbib}
\usepackage{pgf,tikz}
\usepackage{pgfplots}
\usepackage{rotating}
\usepackage{subfig}
\usepackage{url}
\usepackage{hyperref} % cross-refs and links
\hypersetup{colorlinks=true, linkcolor=blue, citecolor=blue, urlcolor=blue}
\usepackage{tikz}
\usetikzlibrary{arrows.meta,positioning,shapes.geometric,fit,backgrounds}

\numberwithin{equation}{section}

\newcommand{\beqn}[1]{\begin{equation}\label{#1}}
\newcommand{\eeqn}{\end{equation}}

\renewcommand{\Pr}{\mathbb{P}}

\usepackage{accents}
\newlength{\dhatheight}

\newtheorem{remark}{Remark}[section]

\usetikzlibrary{arrows}
\tikzstyle{decision} = [diamond, draw, fill=blue!20, text width=4.5em, text badly centered, node distance=3cm, inner sep=0pt]
\tikzstyle{block} = [rectangle, draw, fill=blue!20, text centered, rounded corners, minimum height=4em]
\tikzstyle{line} = [draw, -latex']
\tikzstyle{cloud} = [draw, ellipse, fill=red!20, node distance=3cm, minimum height=2em]
\tikzstyle{cloud2} = [draw, ellipse, fill=green!20, node distance=3cm, minimum height=2em]

\definecolor{aau2}{rgb}{0.0, 0.5, 0.69}
\definecolor{aau3}{rgb}{0.0, 0.53, 0.74}
\definecolor{aau4}{rgb}{0.0, 0.48, 0.65}
\definecolor{aau5}{rgb}{0.0, 0.45, 0.73}
\definecolor{rsap}{RGB}{130, 36, 51}
\definecolor{gsap}{RGB}{112, 164, 137}
\definecolor{ttffqq}{rgb}{0.0, 0.48, 0.65}
\definecolor{ffqqqq}{rgb}{0.0, 0.5, 0.69}
\definecolor{tud}{rgb}{0.43, 0.73, 0.11}
\definecolor{verde}{rgb}{0.33, 0.53, 0.11}
\definecolor{darkgreen}{rgb}{0, 0.6, 0}

\allowdisplaybreaks[1]

\begin{document}

% -----------------------------
% Title & authors
% -----------------------------
\title{On the possibilities and limitations \\ of LLM-based multi-agent systems \\ for multi-criteria decision-making \\ in healthcare operations management}

\author{M. Boresta$^\ddag$ \and T. Giovannelli$^\dag$ \and S. Namirad$^\dag$ \and G. Riccardi$^\S$}
\maketitle

\newcommand{\affiliationUC}{$^\dag$Department of Mechanical, Materials, and Industrial Engineering, University of Cincinnati, Cincinnati, OH}
\newcommand{\affiliationIASI}{$^\ddag$Istituto di Analisi dei Sistemi ed Informatica (IASI), CNR, Rome, Italy}
\newcommand{\affiliationDIAG}{$^\S$Department of Computer, Control, and Management Engineering (DIAG), Sapienza University of Rome, Rome, Italy}
\footnotetext{\hspace{-0.55cm}\affiliationUC}
\footnotetext{\hspace{-0.55cm}\affiliationIASI}
\footnotetext{\hspace{-0.55cm}\affiliationDIAG}
\footnotetext{\hspace{-0.55cm}Emails: marco.boresta@iasi.cnr.it, giovanto@ucmail.uc.edu, namirasr@mail.uc.edu, riccardi@diag.uniroma1.it}

\begin{abstract}
Multi-criteria models for healthcare operations need weights that express how stakeholders trade off patient waiting time, provider workload, and operating cost, but these weights are difficult to quantify. Teams of large language model (LLM) agents, each prompted to represent a stakeholder, have been used to supply them. We propose a framework for assessing the reliability of LLM-generated preferences as inputs to multi-objective optimization and quantifying how their uncertainty affects decisions. We study explicit preferences by asking agents to rate the importance of the criteria in an outpatient appointment scheduling problem, and implicit preferences by inferring trade-offs from the agents' pairwise choices between emergency department configurations. The computational experiments show that the explicit preferences respond to the clinic's operating conditions and differ across roles, but deliberation among the agents can erase these differences, whereas simple averaging preserves them. The pairwise choices used to reveal the agents' implicit preferences depend largely on which option is presented first, and we find limited evidence that a statistical model using contextual information alongside presentation order predicts the agents' choices on new comparison pairs better than a model using presentation order alone. We then formulate robust weighted-sum optimization models to account for preference uncertainty. In most cases, these models select different schedules or configurations from those obtained using average preference weights. These results indicate that LLM-generated preferences need an audit of this kind before they are used as optimization weights in healthcare decision support.
\end{abstract}

%%%%%%%%%%%%%%%%%%%%%%%%%%%%%%%%%%%%%%%%%%%%%%%%%%%%%%%%%%%%%%%%%%%%%%%%%%%%%%%%%%%%%%%%%

\section{Introduction}\label{sec:intro}
%%%%%%%%%%%%%%%%%%%%%%%%%%%%%%%%%%%%%%%%%%%%%%%%%%%%%%%%%%%%%%%%%%%%%%%%%%%%%%%%%%%%%%%%%

Operational decisions in healthcare are inherently multi-criteria. In outpatient clinics, an appointment schedule trades patient waiting time against provider idle time and overtime~\citep{cayirli_outpatient_2003,denton_sequential_2003,gupta_appointment_2008}. In emergency departments (EDs), dedicated pathways for low-acuity patients are widely used to reduce crowding and waiting times, but they require additional staffed opening hours~\citep{boresta_managing_2024}.
Different stakeholders value these criteria differently~\citep{green_queueing_2006,gupta_appointment_2008}: patients primarily care about waiting, providers about workload and overtime, and operations staff about capacity use and cost. Multi-objective optimization characterizes such conflicts through Pareto-optimal solutions, which represent trade-offs among the criteria and provide alternatives for managers to consider~\citep{boresta_managing_2024}. Scalarization methods, such as the weighted-sum method, use preference information to combine multiple criteria into a single objective and guide the search for a preferred solution~\citep{miettinen_nonlinear_2012,ehrgott_multicriteria_2005,abdelaziz_multi-objective_2007,hu_robust_2011,giovannelli_integrated_2023}. In practice, the weights are often chosen informally and then treated as fixed, in healthcare and other scheduling applications~\citep{jiang_integer_2017,maccarrone_simheuristic_2018}, although they may vary across stakeholders and operating conditions and are difficult to state numerically~\citep{gupta_appointment_2008,cayirli_outpatient_2003}.

Large language models (LLMs) offer a new route to these weights. An agent prompted with a stakeholder role can report weights on demand and for each operating condition, whereas eliciting them from patients, providers, and managers is slow and difficult to repeat. Recent work in multi-criteria decision making (MCDM) uses LLM agents to translate qualitative judgments into numerical criterion weights and to reconcile several ``expert'' viewpoints through structured multi-agent procedures~\citep{vahidnia_multi-agent_2025}, uses LLMs as automated evaluators within MCDM pipelines~\citep{wang_one_2025}, or lets an LLM select among multi-objective alternatives directly from stated priorities~\citep{jovine_listen_2026}. Multi-agent LLM systems in which agents play distinct roles are now a general paradigm for planning and coordination~\citep{guo_large_2024}, and a team of agents mirroring the stakeholders of a healthcare system, such as an outpatient clinic or an~ED, is a natural instance of this design. However, it is unclear whether the preferences generated by these agents are reliable enough to be used in multi-objective optimization. LLM outputs vary across repeated queries~\citep{ahn_prompt-reverse_2025,smith_large_2025} and prompt formulations~\citep{ahn_prompt-reverse_2025}; their choices can change when the same two options are presented in reverse order~\citep{zheng_judging_2023,wang_large_2024}; they tend to adapt their answers to the views stated by their users~\citep{sharma_towards_2025}; and agents may drift from their assigned personas toward the model's default positions~\citep{taubenfeld_systematic_2024,chuang_debate_2025}. 
% are not guaranteed to reproduce human behavior~\citep{gao_take_2025};

In high-stakes settings, decision support should be transparent and interpretable~\citep{rudin_stop_2019}, since opaque algorithms can encode hidden biases~\citep{obermeyer_dissecting_2019}. This requirement becomes more pressing as data-driven optimization tools are deployed in healthcare operations~\citep{bertsimas_predictive_2020,bertsimas_hospital-wide_2023}.
Weights generated by LLM agents put this requirement at risk, because they inherit the inconsistencies listed above and no external reference exists against which to check them. We propose a framework for Preference Assessment and Decision
Analysis for LLM-based agents~(LLM-PADA) that audits LLM-generated preferences
before they are used as inputs to multi-objective optimization and quantifies how their
uncertainty affects decisions. The framework considers two
mechanisms. 
The first mechanism is \emph{explicit preference elicitation}, studied in outpatient appointment scheduling. Here, intake, provider, and scheduler agents rate the importance of patient waiting time, provider idle time, and provider overtime on a $1$--$9$ scale. The agents then revise their ratings over three rounds under a numerical protocol, in which they exchange scores, or an argumentative protocol, in which they also exchange justifications and objections. Both are compared with a mechanical protocol that averages the scores without further LLM calls. 
The second mechanism is \emph{implicit preference inference}, studied in an ED in which eligible low-acuity patients can be routed to a minor injuries unit~(MIU). Rather than reporting weights, agents with patient, physician, and operations perspectives choose between pairs of nondominated configurations taken from the discrete-event simulation study of \citet{boresta_managing_2024}, described by their average waiting time and weekly MIU opening hours.

The~LLM-PADA framework has four steps. First, we obtain preferences
through explicit preference elicitation or implicit preference
inference. Second, we assess the reliability of the resulting preferences
by examining how they respond to the decision context and how they are affected by the procedure used to obtain them. Here, reliability refers to internal checks (context sensitivity, role differentiation, and stability under changes in procedure and presentation order), not to agreement with human stakeholder preferences. Third, we construct uncertainty sets for the criterion weights
from the elicited or inferred preferences. Fourth, following the multiexpert multicriteria robust weighted-sum
(McRow) approach of Hu and
Mehrotra~\citep{hu_robust_2011}, we formulate robust weighted-sum optimization models over these sets to quantify
how accounting for preference uncertainty changes decisions.
We demonstrate the~LLM-PADA framework in outpatient appointment scheduling
and emergency department planning using locally executed models
of the Qwen3.5 family~\citep{yang_qwen3_2025}.  Although our computational experiments focus on two healthcare contexts and locally executed models, the framework is general and can be used for other healthcare operations problems and with proprietary~LLMs.

% Note that we make no claim that the outputs of the~LLM agents match the preferences of human stakeholders (i.e., patients, providers, or managers), which is beyond the scope of this paper.

\subsection{Contributions}
This paper makes three contributions:
(i)~we introduce the~LLM-PADA framework, which assesses the reliability of LLM-generated preferences before they are used as objective weights and measures how accounting for their uncertainty changes decisions by comparing mean-weight and robust weighted-sum models; the two models select different decisions in ten of the sixteen outpatient comparisons and in five of the six emergency department comparisons;
(ii)~for explicit preferences in outpatient appointment scheduling, we develop an approach that first elicits importance scores from agents with different stakeholder roles and then compares three ways of revising these scores: two LLM deliberation protocols, in which agents exchange scores alone or scores with arguments, and mechanical averaging; experiments with prompts describing different clinic conditions show that~LLM deliberation can weaken their responses to these conditions and erase or reverse role differences, whereas mechanical averaging retains role differentiation;
(iii)~for implicit preferences in~ED planning, we propose an approach to estimate trade-off weights from agents' choices between pairs of ED configurations described by two outcomes; we use a choice model that accounts for which alternative is presented first and evaluate it on other pairs not used for estimation; we find evidence that outcome differences improve prediction of the agents' choices beyond presentation order for only two of the six combinations of model and role.

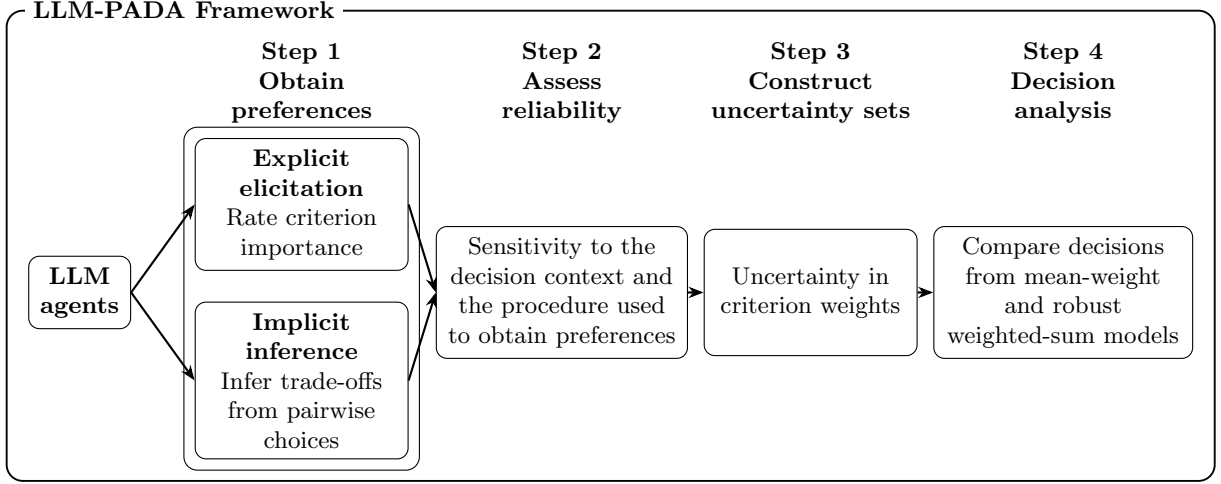
\begin{figure}[htbp]
\centering
\resizebox{\linewidth}{!}{%
\begin{tikzpicture}[
    box/.style={
        draw=black, fill=white, text=black,
        rounded corners=4pt, align=center,
        inner sep=4pt, font=\footnotesize
    },
    src/.style={box},
    mech/.style={box},
    assess/.style={box},
    uncert/.style={box},
    dec/.style={box},
    stephead/.style={
        align=center, text width=2.8cm, text=black,
        font=\footnotesize\bfseries
    },
    arr/.style={-{Stealth[length=2mm]}, draw=black, thick}
]
\node[src, text width=1.05cm, minimum height=0.8cm] (a)
    at (0,0) {\textbf{LLM}\par\textbf{agents}};

\node[mech, text width=2.5cm, minimum height=1.7cm] (s)
    at (2.9,1.15)
    {\textbf{Explicit}\par
     \textbf{elicitation}\par
     Rate criterion\par importance};

\node[mech, text width=2.5cm, minimum height=1.7cm] (r)
    at (2.9,-1.15)
    {\textbf{Implicit}\par
     \textbf{inference}\par
     Infer trade-offs\par from pairwise\par choices};

\node[assess, text width=3.0cm, minimum height=1.65cm] (v)
    at (6.3,0)
    {Sensitivity to the
     decision context and
     the procedure used
     to obtain preferences};

\node[uncert, text width=2.5cm, minimum height=1.65cm] (u)
    at (9.55,0)
    {Uncertainty in
     criterion weights};

\node[dec, text width=3.1cm, minimum height=1.65cm] (d)
    at (12.85,0)
    {Compare decisions\par
     \mbox{from mean-weight}\par
     \mbox{and robust}\par
     weighted-sum models};

\node[stephead, anchor=north] (h1) at (2.9,3.4)
    {Step 1\par Obtain\par preferences};
\node[stephead, anchor=north] (h2) at (6.3,3.4)
    {Step 2\par Assess\par reliability};
\node[stephead, anchor=north] (h3) at (9.55,3.4)
    {Step 3\par Construct\par uncertainty sets};
\node[stephead, anchor=north] (h4) at (12.85,3.4)
    {Step 4\par Decision\par analysis};

\begin{scope}[on background layer]
    % The two preference routes belong to a single step.
    \node[draw=black, rounded corners=5pt,
          inner sep=4pt, fit=(s)(r)] {};
    % The outer frame encloses all four steps and their agent input.
    \node[draw=black, rounded corners=6pt,
          line width=0.7pt, inner sep=8pt,
          fit=(a)(s)(r)(v)(u)(d)(h1)(h2)(h3)(h4)]
          (framework) {};
\end{scope}

\node[anchor=west, fill=white, text=black, inner sep=3pt,
      font=\footnotesize\bfseries]
    at ([xshift=7pt]framework.north west) {LLM-PADA Framework};

\draw[arr] (a.east) -- (s.west);
\draw[arr] (a.east) -- (r.west);
\draw[arr] (s.east) -- (v.west);
\draw[arr] (r.east) -- (v.west);
\draw[arr] (v.east) -- (u.west);
\draw[arr] (u.east) -- (d.west);
\end{tikzpicture}%
}
\caption{Overview of the LLM-PADA framework. Preferences are obtained from LLM agents in two ways, by asking them to rate the criteria in outpatient clinics or by observing their choices between different~ED configurations. In both cases, the resulting weights are uncertain, and we propose robust models to quantify how accounting for preference uncertainty changes decisions.}
\label{fig:paper-overview}
\end{figure}

Table~\ref{tab:literature-gap} in
Section~\ref{sec:related_background} summarizes the gaps in the
literature addressed by these contributions.
Figure~\ref{fig:paper-overview} summarizes the proposed framework.

% Preamble requirements:
% \usepackage{graphicx}
% \usepackage{tikz}
% \usetikzlibrary{arrows.meta,fit,backgrounds}

\subsection{Organization of the paper}
Section~\ref{sec:related_background} reviews related work.
Section~\ref{sec:explicit-preferences} studies explicit
preference elicitation in outpatient clinics.
Section~\ref{sec:implicit_preference} studies implicit preference inference in emergency department planning.
Section~\ref{sec:future_work} concludes and discusses directions for future work.

\section{Background and related work}\label{sec:related_background}
% =========================================================================

\subsection{Multi-objective optimization background}\label{subsec:background}
We first introduce the multi-objective optimization terminology used in the paper. Consider problems of the form $\min_{x\in X} f(x)=(f_1(x),\ldots,f_q(x))$ with $q\ge 2$ objectives, where $X\subseteq\mathbb R^n$ is the feasible set in the decision space, possibly discrete as in our applications, and $f$ maps it into the objective space~$\mathbb R^q$. Given $x,\bar x\in X$, we say that $\bar x$ dominates $x$ (and~$f(\bar x)$ dominates $f(x)$) if $f(\bar x)\le f(x)$ componentwise and $f(\bar x)\ne f(x)$. A solution in~$X$ is Pareto optimal if no $\bar x\in X$ dominates it, and we denote the set of Pareto-optimal (or efficient) solutions by $X_{\mathrm E}$. A point~$f(x)$ is nondominated if~$x\in X_{\mathrm E}$, and the Pareto front is the set of all nondominated points, $\{f(x) ~|~ x\in X_{\mathrm E}\}$. Pareto optimality is thus a property of solutions in the decision space, whereas nondominance is a property of points in the objective space. Any two distinct nondominated points involve a trade-off, since each is better than the other in at least one objective, so choosing among them requires preference information~\citep{ehrgott_multicriteria_2005,miettinen_nonlinear_2012}.

Multi-objective optimization methods compute Pareto-optimal solutions, and hence nondominated points, and they differ in when preference information enters~\citep{ehrgott_multicriteria_2005,miettinen_nonlinear_2012}. A priori methods require preference information before solving, in order to identify a preferred solution. A posteriori methods instead approximate the entire set of nondominated points in the objective space, which captures all trade-offs but can be computationally expensive, and the decision maker then selects a preferred point in a post-processing step. Interactive methods bring in the decision maker's preferences progressively during the optimization, with the goal of identifying a preferred solution.
The weighted-sum method is a standard a priori method that combines the objectives into a single one using weights that reflect the decision maker's preferences, and then solves the resulting single-objective problem~\citep{ehrgott_multicriteria_2005,miettinen_nonlinear_2012}.

\subsection{Related work}

Subsections~\ref{subsub:llm_agents}--\ref{subsub:llm_or_health}
review the related literature. Table~\ref{tab:literature-gap}
positions the paper with respect
to the closest works.

\begin{table}[htbp]
\centering
\caption{Positioning of the paper. \checkmark: addressed; $\circ$: partially addressed; --: not addressed.}
\label{tab:literature-gap}
\small
\resizebox{\textwidth}{!}{%
\begin{tabular}{lcccccc}
\toprule
Work & \makecell{Role-based\\LLM agents} & \makecell{Explicit\\weights} & \makecell{Choice-based\\inference} & \makecell{Reliability\\assessment} & \makecell{Robust decisions\\under weight uncertainty} & \makecell{Healthcare\\operations} \\
\midrule
Vahidnia~\citep{vahidnia_multi-agent_2025} & \checkmark & \checkmark & -- & $\circ$ & -- & -- \\
Jovine et al.~\citep{jovine_listen_2026} & -- & \checkmark & $\circ$ & $\circ$ & -- & -- \\
LLM judges~\citep{zheng_judging_2023,wang_large_2024} & -- & -- & $\circ$ & \checkmark & -- & -- \\
LLM debate~\citep{taubenfeld_systematic_2024,chuang_debate_2025} & \checkmark & -- & -- & \checkmark & -- & -- \\
Hu and Mehrotra~\citep{hu_robust_2011} & -- & $\circ$ & -- & -- & \checkmark & -- \\
PRO~\citep{hu_robust_2015,haskell_ambiguity_2016,vayanos_robust_2021} & -- & -- & $\circ$ & -- & \checkmark & -- \\
Lu et al.~\citep{lu_triageagent_2024} & \checkmark & -- & -- & $\circ$ & -- & $\circ$ \\
Masanneck et al.~\citep{masanneck_triage_2024} & -- & -- & -- & $\circ$ & -- & $\circ$ \\
Boresta et al.~\citep{boresta_managing_2024} & -- & -- & -- & -- & $\circ$ & \checkmark \\
\midrule
This paper & \checkmark & \checkmark & \checkmark & \checkmark & \checkmark & \checkmark \\
\bottomrule
\end{tabular}}
\end{table}

\subsubsection{LLM agents for multi-criteria decision-making}\label{subsub:llm_agents}
Recent surveys document the rapid emergence of LLM-based multi-agent systems as a general paradigm for decision making, highlighting their advantages over single-agent designs in complex, multi-step environments~\citep{guo_large_2024}. A single agent usually works in a loop where the model alternates between reasoning and acting~\citep{yao_react_2023}, and it can call external tools, such as calculators or search engines, through APIs~\citep{schick_toolformer_2023}. Multi-agent systems instead give each agent a narrower role, and they differ in how the agents are organized. In hierarchical designs, a coordinating agent hands subtasks to other agents, as in AutoGen~\citep{wu_autogen_2024} and in layered leader--employee networks~\citep{hao_chatllm_2025}. In peer-to-peer designs, agents talk to each other directly by playing roles~\citep{li_camel_2023}. In protocol-based designs, specialized agents follow prescribed interaction patterns, as in MetaGPT~\citep{hong_metagpt_2024}, ChatDev~\citep{qian_chatdev_2024}, and MDTeamGPT~\citep{chen_mdteamgpt_2026}. A further pattern, particularly relevant under conflicting objectives, is structured deliberation, in which agents represent different viewpoints and, over several rounds, critique and revise each other's positions~\citep{TGiovannelli_GDKent_2026}. For such systems, performance should be assessed not only by end-to-end task success but also by properties such as consistency across repeated runs~\citep{yehudai_survey_2026}. For deliberative systems, the agreement induced by the protocol is a further property of interest~\citep{taubenfeld_systematic_2024,chuang_debate_2025}. \textit{We use this deliberative design for explicit preferences.} The agents revise their scores for three rounds, and no single agent makes the final decision.

A growing literature uses LLM agents to produce the preference information that multi-criteria decision models require. In multi-criteria decision making (MCDM), Vahidnia~\citep{vahidnia_multi-agent_2025}, whose procedure inspired our explicit elicitation, obtains initial criterion weights from multiple LLM instances acting as experts and then enforces agreement through an aggregation-and-convergence procedure that removes outliers and iteratively moves the group weights toward the global mean until the well-known Kendall's coefficient of concordance stabilizes. This convergence step is applied mechanically rather than through LLM deliberation, and is therefore close to the mechanical averaging baseline of our explicit study. Wang et al.~\citep{wang_one_2025} propose a general LLM-based MCDM framework in which LLMs replace human experts in rating the alternatives on each criterion, while the criterion weights are obtained with the analytic hierarchy process. They show that off-the-shelf models reach only moderate accuracy against expert judgments, whereas fine-tuning open-source models on the task brings them close to expert level. Jovine et al.~\citep{jovine_listen_2026} introduce LISTEN, an LLM framework for selecting one item from a large set of alternatives described by multiple attributes. LISTEN-U asks an LLM to translate a natural-language description of priorities into a linear utility function and to revise its coefficients iteratively, and is therefore similar to our explicit-preference protocol, but it does not involve deliberation among agents representing different stakeholder roles. LISTEN-T samples batches of alternatives, asks the LLM to select a winner from each batch, and conducts a final tournament among the winners. Both LISTEN-T and our implicit-preference protocol obtain information through choices, but LISTEN-T uses the choices directly to select a candidate, whereas we treat them as observations, estimate a choice model that separates attribute effects from presentation effects, quantify the uncertainty of the inferred trade-off, and evaluate the resulting decisions through mean-weight and robust optimization. In all these works, the LLM output is treated as the preference to be used. \textit{We ask instead whether it is reliable enough to be used, and what happens to the decisions when it is not.}

\subsubsection{Reliability of LLM-generated preferences}
Studies of LLM behavior raise concerns about using agents as reliable sources of stakeholder preferences. Gao et al.~\citep{gao_take_2025} warn against using LLMs as surrogates for human participants, showing that in a strategic reasoning game their behavior can differ markedly from that of humans. Additionally, LLM answers can be inconsistent beyond the variability induced by sampling and paraphrasing. For example, Ahn and Yin~\citep{ahn_prompt-reverse_2025} show that a model may give conflicting answers when asked which candidates are correct and which are incorrect, and an informal test in which the same questions were posed twice to several commercial chatbots produced contradictory answers~\citep{smith_large_2025}. In medical tasks, the consistency of LLM outputs depends on prompt design~\citep{wang_prompt_2024}, and suboptimal repeatability has been reported even for a structured disaster-triage protocol~\citep{franc_accuracy_2024}, although accuracy on clinical prediction tasks has been found to be largely insensitive to sampling temperature~\citep{patel_exploring_2024}. Interaction among agents raises a further concern. Multi-agent debate, in which agents critique and revise one another's answers over several rounds to reach a common answer, has been proposed to improve factuality and reasoning~\citep{du_improving_2023}. Agreement, however, is not necessarily informative, as the following examples show: in simulated political debates, role-playing agents drift from their assigned personas toward the model's own default positions~\citep{taubenfeld_systematic_2024}; models trained with human feedback tend to adapt their answers to the views stated by their users~\citep{sharma_towards_2025}; and in simulated debates groups of LLM agents converge more strongly than groups of humans~\citep{chuang_debate_2025}. 
These studies concern factual questions or opinions on social topics. In a team of stakeholder agents, the same tendencies could erase the very differences in priorities that justify having a team. \textit{Our explicit-preference study therefore measures agreement separately from whether the agents keep their contextual and role-specific priorities}, and it compares LLM deliberation with a mechanical averaging rule that involves no LLM calls.

\subsubsection{Preference inference from pairwise choices}

Our implicit-preference study builds on classical models of pairwise choice. The Bradley--Terry model assigns a latent score to each alternative and models the probability that one is preferred to the other as a logistic function of the difference between their scores~\citep{bradley_rank_1952}. Random-utility models with logit choice probabilities extend this by expressing the probabilities in terms of the alternatives' attributes~\citep{mcfadden_measurement_1974}. The reward models used to align LLMs with human feedback have the same structure and are fitted to pairwise human comparisons~\citep{zeng_aligning_2025}. We use the structure the other way around: the LLM agent is the respondent, the alternatives are simulated~ED operating configurations described by two outcomes, and the ratio of the fitted outcome coefficients gives the weight in a weighted-sum objective. Studies of LLMs as judges of pairwise comparisons find that their verdicts depend on the order in which the alternatives are presented~\citep{zheng_judging_2023,wang_large_2024}. \textit{Our choice model therefore includes a presentation term}, and we assess its predictive value against a benchmark that uses presentation order alone.

\subsubsection{Robust optimization under preference uncertainty}

When the weights used to represent a decision maker's preferences
and combine objectives into a single objective are uncertain, robust optimization~\citep{ben-tal_robust_2009} protects against the worst case over a set of plausible weights or utilities. Minmax robustness for multi-objective problems with uncertain parameters is formalized by Ehrgott et al.~\citep{ehrgott_minmax_2014}, and Hu and Mehrotra~\citep{hu_robust_2011} study robust and stochastically weighted multiobjective models in which the weights belong to an uncertainty set. Our robust problems are deterministic instances of their multi-expert multi-criteria robust weighted sum (McRow) approach. A related literature, known as preference robust optimization~(PRO), constructs ambiguity sets of utility functions from partial information elicited from the decision maker. Hu and Mehrotra~\citep{hu_robust_2015} formulate a maximin model over a set of risk-averse utilities, Wu et al.~\citep{wu_bi-attribute_2025} approximate bi-attribute utilities by piecewise-linear functions with error bounds, and Guo et al.~\citep{guo_utility_2024} build ambiguity sets from moment-type conditions. Hu et al.~\citep{hu_distributional_2025} model the decision maker's utility as random with an ambiguous distribution, and Haskell et al.~\citep{haskell_ambiguity_2016} consider ambiguity in the risk preferences and in the underlying probability distribution simultaneously, whereas Natarajan et al.~\citep{natarajan_tractable_2010} assume a known piecewise-linear concave utility and robustify only against ambiguity in the return distribution. Vayanos et al.~\citep{vayanos_robust_2021} combine a robust model with linear utilities in a polyhedral uncertainty set and an elicitation process that refines the set from the answers to pairwise comparison queries. We also use linear utility functions, but represent uncertainty in their weights using an empirical percentile box in the outpatient study and a bootstrap interval in the emergency department study. Whereas Vayanos et al.~\citep{vayanos_robust_2021} select comparison queries to maximize the worst-case utility, \textit{we study how variability and uncertainty in LLM-generated weights affect operational decisions}.

\subsubsection{LLMs in operations research and healthcare}\label{subsub:llm_or_health}
Our work relates to two applied areas. In operations research, LLMs are increasingly embedded in optimization pipelines with closed-loop evaluation, where objective feedback is used to filter or refine candidate solutions. FunSearch~\citep{romera-paredes_mathematical_2024} and HeuriGym~\citep{chen_heurigym_2026} evaluate LLM-generated programs and heuristics by executing them, Çetinkaya et al.~\citep{cetinkaya_discovering_2025} derive scheduling heuristics and validate them against exact methods, and Xiao et al.~\citep{xiao_chain--experts_2024} coordinate domain-specialized agents that formulate and solve operations research problems. In healthcare, LLMs and multi-agent LLM systems have been evaluated for clinical triage~\citep{lu_triageagent_2024,masanneck_triage_2024}, for broader triage, referral, and diagnosis workflows~\citep{gaber_evaluating_2025}, for medical consultation~\citep{ren_healthcare_2025}, for transitions from primary to specialist care~\citep{tao_llm_2026}, and for summarizing patient dialogues~\citep{ferreira_comprehensive_2025}. In both areas, the LLM output can be checked against a solver, a simulator, a clinical reference, or human evaluation. \textit{The objective weights in this paper have no such reference}, so we focus on diagnosing how reliable they are and on measuring their effect on decisions.

%%%%%%%%%%%%%%%%%%%%%%%%%%%%%%%%%%%%%%%%%%%%%%%%%%%%%%%%%%
% SECTION 3
%%%%%%%%%%%%%%%%%%%%%%%%%%%%%%%%%%%%%%%%%%%%%%%%%%%%%%%%%%
\section{Explicit preference elicitation in outpatient scheduling}
\label{sec:explicit-preferences}

Outpatient appointment scheduling involves trade-offs among patient
waiting time, provider idle time, and provider overtime~\citep{cayirli_outpatient_2003,denton_sequential_2003,gupta_appointment_2008}. Patients, providers, and operations staff weigh these criteria differently, and the weights may vary with clinic conditions. This section applies the~LLM-PADA framework to study whether a team of LLM agents, each prompted with one stakeholder role, can supply such weights explicitly by reporting a numerical importance score for each criterion, and whether these scores are reliable enough to serve as coefficients of a scheduling objective. We use \emph{preference} for the importance scores that an agent reports under the experimental prompts; whether these scores match the preferences of human patients, providers, or schedulers is outside the scope of the paper.

%%%%%%%%%%%%%%%%%%%%%%%%%%%%%%%%%%%%%%%%%%%%%%%%%%%%%%%%%%
\subsection{Problem setting}
\label{sec:explicit-setting}

We consider an outpatient clinic session in which one provider must schedule patient appointment requests. The session has length $H$ minutes and uses appointment slots of length $\tau>0$ minutes, indexed by $\mathcal T=\{0,\ldots,m-1\}$.
Let $\mathcal P=\{1,\ldots,n\}$ denote the set of $n$ patient requests, listed in the order in which they must be scheduled. Each request $p\in\mathcal P$ requires $d_p$ minutes of service and cannot be scheduled before slot $a_p\in\mathcal T$. The scheduling decision assigns one appointment slot $x_p\in\mathcal T$ to each request $p$, and we collect these assignments in the vector $x=(x_p)_{p\in\mathcal P}$. We denote by~$X$ the set of feasible schedules:
\begin{equation}
 X=\Bigl\{x\in\mathcal T^{n}:\ x_p\ge a_p,\ \ \forall\, p\in\mathcal P;\quad
 x_{p+1}\ge x_p+1,\ \ \forall\, p\in\mathcal P\setminus\{n\}\Bigr\}.
 \label{eq:explicit-feasible-schedules}
\end{equation}
The first family of constraints ensures that each patient is scheduled no earlier than its earliest acceptable slot. The second preserves the prescribed order of the requests and prevents more than one request from being assigned to the same slot. Service durations are deterministic and known, the request order is fixed, and overbooking is not allowed. These choices remove sources of variation other than the weights, whose effect we study. Stochastic service durations and overbooking are central to the appointment scheduling literature~\citep{cayirli_outpatient_2003,gupta_appointment_2008}, but they are not needed to isolate this effect.

For a feasible schedule $x\in X$, patient $p$ is assigned to appointment time $t_p=\tau\,x_p$. Because a service may extend beyond the next patient's appointment time, the actual service start time can differ from the scheduled appointment time. The first patient starts service at $b_1=t_1$, while each subsequent patient starts at
$b_p=\max\{t_p,b_{p-1}+d_{p-1}\}$ for $p\in\mathcal P\setminus\{1\}$. Thus, patient $p$ waits $q_p=b_p-t_p$ minutes. The provider is idle for $g_1=t_1$ minutes before the first patient and for $g_p=b_p-b_{p-1}-d_{p-1}$ minutes before patient $p\in\mathcal P\setminus\{1\}$. If the last service extends beyond the session length $H$, the resulting overtime is $O=\max\{0,b_n+d_n-H\}$. Let $\nu_w,\nu_u,\nu_o>0$ denote fixed reference scales in minutes. We define the following three normalized objective functions
\begin{equation}
 f_w(x)=\frac{1}{\nu_w}\sum_{p\in\mathcal P}q_p,
 \qquad
 f_u(x)=\frac{1}{\nu_u}\sum_{p\in\mathcal P}g_p,
 \qquad
 f_o(x)=\frac{O}{\nu_o},
 \label{eq:explicit-scheduling-objectives}
\end{equation}
and write $f(x)=(f_w(x),f_u(x),f_o(x))^\top$. Idle time includes the initial gap before the first service and the gaps between consecutive services, but not the unused time after the last service. The resulting multi-objective appointment scheduling problem is
\begin{equation}
 \min_{x\in X}\; f(x).
 \label{eq:explicit-multiobjective}
\end{equation}

\subsection{Research questions}
\label{sec:research_question}
We scalarize problem~\eqref{eq:explicit-multiobjective} with the weighted-sum method~\citep{miettinen_nonlinear_2012,ehrgott_multicriteria_2005}, whose criterion weights represent trade-offs among the objectives, and we obtain these weights from a team of LLM agents representing different stakeholder perspectives. In Step~1 of LLM-PADA, each agent first reports numerical importance scores for the criteria (\emph{elicitation}); the agents then observe one another's scores and revise their own over several rounds (\emph{deliberation}). We address four research questions.
\begin{itemize}
 \item[(Q1)] \emph{Contextual sensitivity:} Do an agent's scores change as expected when the clinic's operational state varies?
\item[(Q2)] \emph{Role differentiation:} After deliberation,
does each agent still assign greater relative importance to
the criterion associated with its role than the other agents do?
\item[(Q3)] \emph{Inter-agent dispersion:} How much do the agents'
normalized weight vectors differ, and how do these differences
change across deliberation rounds?
\item[(Q4)] \emph{Scheduling consequences:} How does accounting for uncertainty in the agents' weights change the selected appointment schedule compared with using mean weights?
\end{itemize}

Low inter-agent dispersion~(Q3) means that the agents agree, but agreeing agents may still fail to respond to clinic conditions~(Q1) or lose their role-specific priorities~(Q2). We therefore assess (Q1)--(Q3) separately before studying the scheduling consequences of the resulting weights~(Q4).

%%%%%%%%%%%%%%%%%%%%%%%%%%%%%%%%%%%%%%%%%%%%%%%%%%%%%%%%%%
\subsection{Methodology}
\label{sec:explicit-methodology}

We now describe how the four steps of LLM-PADA are implemented for the outpatient appointment scheduling problem. Step~1 produces the initial and revised preference weights. Step~2 assesses contextual sensitivity, role differentiation, and inter-agent dispersion to answer (Q1), (Q2), and (Q3), respectively. Step~3 constructs uncertainty sets for the weights, and Step~4 compares mean-weight and robust scheduling decisions to answer (Q4).

%%%%%%%%%%%%%%%%%%%%%%%%%%%%%%%%%%%%%%%%%%%%%%%%%%%%%%%%%%
\subsubsection{Preference generation (Step 1)}
\label{sec:explicit-generation}

Step~1 consists of the elicitation and deliberation stages described below. For deliberation, we compare three protocols: a \emph{numerical} protocol, in which the agents exchange score vectors only; an \emph{argumentative} protocol, in which they also exchange short written justifications and objections; and a \emph{mechanical} protocol, in which a fixed rule moves each agent's scores toward the group mean without LLM calls.

\label{sec:explicit-design}
Let $\mathcal E=\{\text{Intake},\text{Provider},\text{Scheduler}\}$ denote the set of agents and $\mathcal C=\{w,u,o\}$ the set of criteria, where $w$, $u$, and $o$ denote patient waiting time, provider idle time, and provider overtime, respectively. In the elicitation stage, each agent receives a prompt specifying its stakeholder role.
The role descriptions in these prompts emphasize timely patient
access for the intake agent, a sustainable clinical workload
for the provider agent, and efficient use of clinic capacity
for the scheduler agent. Every agent scores all three criteria with the same definitions and the same scale. The role descriptions are fixed throughout elicitation and deliberation.

Let $\mathcal I$ index the base clinic descriptions. Each description lists nine fields: current time, session end, patients seen, patients remaining, open slots, mean waiting time, provider idle time, projected overtime, and cancellations. For each description~$i\in\mathcal I$, we create four versions,
one for each condition~$S\in\mathcal S=\{\mathrm{B},\mathrm{HW},\mathrm{HO},\mathrm{CF}\}$.
B (balanced) uses the base description without modification and
serves as the reference. HW (high waiting) modifies the base
description by increasing the displayed mean waiting time.
HO (high overtime) modifies it by moving the displayed current
time later and increasing the projected overtime.
CF (context-free) contains no clinic description, as the prompt
includes only the role description and the criteria.
The numerical modifications and repetition counts are specified
in Subsection~\ref{sec:explicit-experiment-setup}.

For every $i\in\mathcal I$, conditions B, HW, and HO are built from the same base description and differ only in the fields that we manipulate. We call this a \emph{paired design}, as the difference between an agent's scores under HW (or HO) and under B for the same $i$ isolates the response to the manipulation. For example, the difference between HW and B is attributable to the displayed mean waiting time alone, while the difference between HO and B reflects the joint change in current time and projected overtime.

Now, let~$\mathcal J$ index repeated responses. For each agent $e\in\mathcal E$, base description $i\in\mathcal I$, condition $S\in\mathcal S$, and repetition $j\in\mathcal J$, the agent returns one \emph{score vector} $s_e^{(i,S,j)}=\bigl(s_{e,w}^{(i,S,j)},s_{e,u}^{(i,S,j)},s_{e,o}^{(i,S,j)}\bigr)$, whose three components are integers from $1$ to $9$ elicited jointly in one response. The total number of elicitation responses is $|\mathcal E|\,|\mathcal I|\,|\mathcal S|\,|\mathcal J|$.

Fix $(e,c,i,S)\in\mathcal E\times\mathcal C\times\mathcal I\times\mathcal S$. Let $Q_1$ and $Q_3$ be the first and third empirical quartiles of the repeated scores $s_{e,c}^{(i,S,j)}$, $j\in\mathcal J$. We retain the repetitions in~$
 \mathcal J_{e,c}^{(i,S)}=\bigl\{j\in\mathcal J:\ Q_1-1.5\,(Q_3-Q_1)\le s_{e,c}^{(i,S,j)}\le Q_3+1.5\,(Q_3-Q_1)\bigr\}$,
which is the interquartile-range rule of~\citet{tukey_exploratory_1977} for removing occasional aberrant responses. The rule is applied to each criterion separately, even though the three scores of a response are elicited jointly. Hence, $\mathcal J_{e,c}^{(i,S)}$ can differ across criteria. The retained mean is $\bar s_{e,c}^{(i,S)}=|\mathcal J_{e,c}^{(i,S)}|^{-1}\sum_{j\in\mathcal J_{e,c}^{(i,S)}}s_{e,c}^{(i,S,j)}$, and the \emph{normalized weight} of agent $e$ for criterion $c$ is
\begin{equation}
 \lambda_{e,c}^{(i,S)}=
 \frac{\bar{s}_{e,c}^{(i,S)}}{\sum_{c'\in\mathcal C}\bar{s}_{e,c'}^{(i,S)}}.
 \label{eq:explicit-normalization}
\end{equation}
Thus, $\lambda_e^{(i,S)}=(\lambda_{e,w}^{(i,S)},\lambda_{e,u}^{(i,S)},\lambda_{e,o}^{(i,S)})$ lies in the simplex $\Delta=\{\lambda\in\mathbb R_+^3:\sum_{c\in\mathcal C}\lambda_c=1\}$.
In Subsection~\ref{sec:explicit-optimization}, we use the normalized weights as coefficients of the criteria in a weighted-sum scheduling objective. Rating-based weighting is an established approach to criteria weighting in multi-criteria decision analysis (MCDA)~\citep{marsh2016multiple} and has also been employed for LLM-based weight elicitation~\citep{vahidnia_multi-agent_2025}. 
% Because the prompts show neither the objective function nor candidate schedules, this use assumes that the importance an agent reports for a criterion carries over to the corresponding normalized objective in~\eqref{eq:explicit-scheduling-objectives}. 
% In Subsection~\ref{sec:explicit-optimization}, we examine how the selected schedule changes when the weights vary over the range observed across agents, conditions, and chains.

For each base description $i$ and condition $S$, the \emph{initial score matrix} collects the retained means,
\begin{equation}
 M^{(i,S,0)}=
 \begin{bmatrix}
 \bar s_{\mathrm{Intake},w}^{(i,S)} & \bar s_{\mathrm{Provider},w}^{(i,S)} & \bar s_{\mathrm{Scheduler},w}^{(i,S)}\\
 \bar s_{\mathrm{Intake},u}^{(i,S)} & \bar s_{\mathrm{Provider},u}^{(i,S)} & \bar s_{\mathrm{Scheduler},u}^{(i,S)}\\
 \bar s_{\mathrm{Intake},o}^{(i,S)} & \bar s_{\mathrm{Provider},o}^{(i,S)} & \bar s_{\mathrm{Scheduler},o}^{(i,S)}
 \end{bmatrix}.
 \label{eq:explicit-initial-matrix}
\end{equation}
Column $e$ of $M^{(i,S,0)}$ is used as the initial score vector $s_e^{(0)}=\bar s_e^{(i,S)}$ of agent $e$ for every deliberation protocol, with normalized weight vector $\lambda_e^{(0)}=\lambda_e^{(i,S)}$.

\label{sec:explicit-deliberation}
The deliberation stage fixes $(i,S)$, which we omit from the notation.
% , and starts every protocol from the same initial score vectors $s_e^{(0)}$, $e\in\mathcal E$.
The protocols run for $R$ synchronous rounds and differ only in their update rules. In round $r\in\{1,\ldots,R\}$, agent $e$ receives its unchanged role description, the clinic description (for B, HW, and HO), and the three score vectors produced in round $r-1$. All agents update simultaneously, so no agent observes another agent's round-$r$ revision. We write $s_e^{(r)}$ for the score vector after round $r$ and $\lambda_e^{(r)}$ for its normalization by~\eqref{eq:explicit-normalization}.

In the \emph{numerical protocol}, the agent returns its revised score vector only, with integer components from $1$ to $9$. In the \emph{argumentative protocol}, in addition to the revised score vector, the agent returns (i)~a one-sentence justification, (ii)~an objection addressed to one other agent, (iii)~the identifiers of the clinic-state fields that it used, and (iv)~a reference to a preceding argument, if it relies on one. From round~2 onward, the argumentative prompt also contains the three argument records of the previous round, each with its validation status. A response is \emph{valid} if its scores are integers from $1$ to $9$ and, in the argumentative protocol, if it contains all required fields, addresses its objection to another agent, and references only previous-round arguments by other agents. An invalid response is retried; once the permitted retries are exhausted, the failure is recorded and the agent's previous-round score vector is carried forward.

In the \emph{mechanical protocol}, which does not use LLM agents, scores are updated as follows
\begin{equation}
 s_e^{(r+1)}=(1-\alpha)\,s_e^{(r)}+\alpha\,\bar{s}^{(r)},
 \qquad \bar{s}^{(r)}=\frac{1}{|\mathcal E|}\sum_{e'\in\mathcal E}s_{e'}^{(r)},
 \qquad 0<\alpha<1,
 \label{eq:mechanical-update}
\end{equation}
where $\bar s^{(r)}$ is the mean score vector across agents, and $\alpha$ specifies the degree of averaging applied in each round. The rule preserves $\bar s^{(r)}$ and moves each agent's score vector toward this mean. It shows how much agreement, and how much attenuation of contextual and role-specific differences, arises from averaging alone, with no agent exercising any judgment. The protocol is close to the convergence step of~\citet{vahidnia_multi-agent_2025}, which also moves the agents' weights toward the group mean by a fixed rule rather than through deliberation among the agents.

The numerical and argumentative protocols, hereafter referred to as the LLM protocols, use the same scoring scale and also provide the agents with the same clinic descriptions.
For each pair $(i,S)$, we generate $K$ independent deliberation trajectories under each LLM protocol, using different random seeds. In each of the $R$ rounds, every agent produces one response. Therefore, each LLM protocol generates~$|\mathcal E|\cdot R\cdot K\cdot|\mathcal I|\cdot|\mathcal S|$ agent responses in total.
The mechanical protocol is deterministic and therefore produces a single $R$-round trajectory from the initial score vectors associated with each~$(i,S)$. We note that although these initial scores depend on $(i,S)$, the mechanical updates themselves do not use the clinic-state description.

%%%%%%%%%%%%%%%%%%%%%%%%%%%%%%%%%%%%%%%%%%%%%%%%%%%%%%%%%%
\subsubsection{Reliability assessment (Step 2)}
\label{sec:explicit-reliability}

For $S\in\{\mathrm{HW},\mathrm{HO}\}$, we measure the response of agent $e$ to the manipulation by the paired score difference
\begin{equation}
 \delta_{e,c}^{(i)}(S)=
 \bar{s}_{e,c}^{(i,S)}-\bar{s}_{e,c}^{(i,\mathrm{B})},
 \qquad \forall\, e\in\mathcal E,\; c\in\mathcal C,\; i\in\mathcal I.
 \label{eq:explicit-context-shift}
\end{equation}
We refer to the \emph{focal criterion} of a condition as the criterion whose displayed information was changed. We expect a positive difference for the focal criterion.

Let $c(e)$ denote the criterion associated with the role of agent $e$: $c(\text{Intake})=w$, $c(\text{Provider})=o$, and $c(\text{Scheduler})=u$, following the references to timely access, to work beyond the session end, and to unused appointment gaps in the three role descriptions (Appendix~\ref{app:explicit-prompts}). The \emph{role differentiation score} of agent $e$ is
\begin{equation}
 \mathrm{RD}_{e}^{(i,S)}=
 \lambda_{e,c(e)}^{(i,S)}-
 \frac{1}{2}\sum_{e'\in\mathcal E\setminus\{e\}}\lambda_{e',c(e)}^{(i,S)}.
 \label{eq:role-persistence}
\end{equation}
A positive score means that agent $e$ assigns more relative importance to its role-associated criterion than the average of the other two agents does. The prompts ask each agent to remain faithful to its role but do not require $\mathrm{RD}_{e}^{(i,S)}>0$. We compute $\delta_{e,c}^{(i)}(S)$ and $\mathrm{RD}_{e}^{(i,S)}$ from the initial scores and again after round~$R$ of each protocol, with $s_e^{(R)}$ and $\lambda_e^{(R)}$ in place of $\bar s_e^{(i,S)}$ and $\lambda_e^{(i,S)}$.

In every round $r$, the \emph{inter-agent dispersion} is
\begin{equation}
 D^{(r)}=\Biggl(\frac{1}{|\mathcal E|\,|\mathcal C|}\sum_{e\in\mathcal E}\sum_{c\in\mathcal C}\bigl(\lambda_{e,c}^{(r)}-\bar\lambda_{c}^{(r)}\bigr)^2\Biggr)^{1/2},
 \qquad \bar\lambda_{c}^{(r)}=\frac{1}{|\mathcal E|}\sum_{e\in\mathcal E}\lambda_{e,c}^{(r)},
 \label{eq:explicit-dispersion}
\end{equation}
that is, the root-mean-square deviation of the normalized weights from their means across agents. A decline in $D^{(r)}$ indicates greater agreement among agents. For the LLM protocols, all three measures are averaged over the $K$ trajectories of each $(i,S)$ before the bootstrap analysis. 
% Kendall's coefficient of concordance $W$~\citep{kendall_problem_1939}, which measures agreement among the agents' rankings of the criteria, is reported as a secondary diagnostic in Appendix~\ref{app:explicit-results}.

We assess variation in the diagnostic measures by resampling base clinic descriptions with replacement and recomputing the measures. Results for paired conditions are resampled together. The bootstrap settings are specified in Subsection~\ref{sec:explicit-experiment-setup}.

%%%%%%%%%%%%%%%%%%%%%%%%%%%%%%%%%%%%%%%%%%%%%%%%%%%%%%%%%%
\subsubsection{Weight uncertainty construction (Step 3)}
\label{sec:explicit-uncertainty}

Let $\Pi=\{\text{pre},\text{num},\text{arg},\text{mech}\}$ denote four weight sources. Specifically, for $\pi=\text{pre}$, the weights are the normalized initial score vectors $\lambda_e^{(i,S)}$ defined in~\eqref{eq:explicit-normalization}, obtained from the retained mean elicitation scores before any deliberation. For $\pi\in\{\text{num},\text{arg},\text{mech}\}$, the weights are obtained by normalizing each agent's score vector after round $R$ of the corresponding deliberation protocol. For each $(S,\pi)\in\mathcal S\times\Pi$, let $\Lambda^{S,\pi}$ denote the sample of normalized weight vectors associated with condition $S$ and source $\pi$, pooled across agents and base descriptions $i\in\mathcal I$, and, for the numerical and argumentative protocols, across the $K$ deliberation trajectories.
Thus, the pre-deliberation and mechanical samples each contain $|\mathcal E|\,|\mathcal I|$ vectors, whereas the numerical and argumentative samples each contain $|\mathcal E|\,|\mathcal I|\,K$ vectors. This gives $|\mathcal S|\,|\Pi|$ samples.

For each criterion $c\in\mathcal C$, let $L_c^{S,\pi}$ and $U_c^{S,\pi}$ be the 5th and 95th empirical percentiles of the $c$-th components of the vectors in $\Lambda^{S,\pi}$, and let $\bar\lambda^{S,\pi}$ be the arithmetic mean of these vectors. We define the corresponding weight set as
\begin{equation}
 \mathcal U_{S,\pi}=\Bigl\{\lambda\in\Delta:
 L_c^{S,\pi}\le\lambda_c\le U_c^{S,\pi},\ \ \forall\, c\in\mathcal C\Bigr\}.
 \label{eq:explicit-uncertainty-set}
\end{equation}

We call $\mathcal U_{S,\pi}$ an empirical percentile-box weight uncertainty set. It is a nonempty polyhedron by construction and reflects variation across stakeholder agents, base clinic descriptions, and, for the LLM protocols, deliberation trajectories. Because the percentile bounds are computed separately for each criterion, $\mathcal U_{S,\pi}$ can contain weight vectors whose components were never observed together in~$\Lambda^{S,\pi}$. We interpret $\mathcal U_{S,\pi}$ as an empirical set of plausible weights and use it in the robust problem~(RB) in Section~\eqref{sec:explicit-optimization} below to guard against every weight vector within the observed componentwise ranges, including combinations that no agent reported.

%%%%%%%%%%%%%%%%%%%%%%%%%%%%%%%%%%%%%%%%%%%%%%%%%%%%%%%%%%
\subsubsection{Scheduling decision analysis (Step 4)}
\label{sec:explicit-optimization}

In Step~4, we compare schedules obtained using mean weights with those obtained by robust optimization to assess how accounting for weight uncertainty affects the scheduling decision.
For each $(S,\pi)\in\mathcal S\times\Pi$, we solve
\begin{equation}
 \text{(MW)}\quad\min_{x\in X}\ (\bar\lambda^{S,\pi})^\top f(x),
 \qquad\qquad
 \text{(RB)}\quad\min_{x\in X}\ \max_{\lambda\in\mathcal U_{S,\pi}}\ \lambda^\top f(x),
 \label{eq:explicit-problems}
\end{equation}
where~$f(x)$ was defined in~\eqref{eq:explicit-scheduling-objectives}. The mean-weight problem (MW) uses the arithmetic mean of the
normalized weight vectors in the corresponding sample. Following
Hu and Mehrotra~\citep{hu_robust_2011}, the robust problem
(RB) is a deterministic multiexpert multicriteria robust
weighted-sum (McRow) model. For each feasible schedule $x$,
the inner maximization identifies the weight vector in
$\mathcal U_{S,\pi}$ that gives the largest weighted objective.
The outer minimization then selects the schedule with the
smallest such value. The service durations
and scheduling constraints are fixed, so robustness concerns
uncertainty in the criterion weights.

Let $\mathcal V_{S,\pi}=\operatorname{vert}(\mathcal U_{S,\pi})$ denote the finite set of vertices of the weight region. Because $\lambda^\top f(x)$ is linear in $\lambda$, the inner maximization in (RB) is attained at a vertex. The McRow problem therefore admits the following
mixed-integer linear programming formulation:
\begin{equation}
 \begin{aligned}
 \min_{x\in X,\,b,\,O,\,\zeta}\quad & \zeta\\
 \text{s.t.}\quad
 & \zeta\ge(\lambda^{v})^\top f(x,b,O), && \forall\, v\in\mathcal V_{S,\pi},\\
 & b_p\ge \tau\,x_p, && \forall\, p\in\mathcal P,\\
 & b_p\ge b_{p-1}+d_{p-1}, && \forall\, p\in\mathcal P\setminus\{1\},\\
 & O\ge b_n+d_n-H,\quad O\ge 0,
 \end{aligned}
 \label{eq:explicit-milp}
\end{equation}
where $\lambda^{v}$ is the weight vector associated with vertex $v$. The criterion functions are
$f_w(x,b,O)=\frac{1}{\nu_w}\sum_{p\in\mathcal P}(b_p-\tau\,x_p)$,
$f_u(x,b,O)=\frac{1}{\nu_u}\bigl(b_n-\sum_{p\in\mathcal P\setminus\{n\}}d_p\bigr)$, and
$f_o(x,b,O)=O/\nu_o$. Problem~\eqref{eq:explicit-milp} is a mixed-integer linear program (MILP), with integer variables $x$ and continuous variables $b$, $O$, and $\zeta$. Problem (MW) uses the same constraints but replaces $\zeta$ and the vertex constraints with the objective $(\bar\lambda^{S,\pi})^\top f(x,b,O)$.

To answer (Q4), we compare the schedules selected by (MW) and (RB) in patient waiting time, provider idle time, and overtime, and evaluate their weighted objectives at the mean weights and over $\mathcal U_{S,\pi}$. The reliability results for (Q1)--(Q3) guide the interpretation of these decision comparisons.

%%%%%%%%%%%%%%%%%%%%%%%%%%%%%%%%%%%%%%%%%%%%%%%%%%%%%%%%%%
\subsection{Numerical experiments and results}
\label{sec:explicit-results}

%%%%%%%%%%%%%%%%%%%%%%%%%%%%%%%%%%%%%%%%%%%%%%%%%%%%%%%%%%
\subsubsection{Experimental setup}
\label{sec:explicit-experiment-setup}

We implement the design of Subsection~\ref{sec:explicit-methodology} with the locally executed model Qwen3.5-4B for all three agents. Appendix~\ref{app:explicit-prompts} provides the prompt templates and validation rules.

We use ten synthetic base clinic descriptions,
$\mathcal I=\{1,\ldots,10\}$, and collect twenty elicitation responses,
$\mathcal J=\{1,\ldots,20\}$, for each agent, description, and condition.
Following the design in Subsection~\ref{sec:explicit-design}, HW adds
25--40 minutes to the displayed mean waiting time, while HO moves the
displayed current time two hours later, capped at 30 minutes before
the session end, and adds 35--55 minutes to projected overtime.
For CF, the same role-and-criteria prompt is used for all ten values
of $i$, without a clinic description. Elicitation uses temperature
$0.7$ and produces $2{,}400$ responses in total.

Deliberation proceeds for $R=3$ rounds from the same initial score
matrix under all three protocols. For each $(i,S)$, we generate
$K=5$ independent trajectories under each LLM protocol, using
different random seeds and temperature $0.2$. This gives $1{,}800$
agent responses per LLM protocol, excluding retries. We use the
lower temperature to reduce sampling variability across rounds so
that score changes more clearly reflect the information exchanged
during deliberation. The mechanical protocol uses $\alpha=0.35$, so each round moves every agent's score 35\% of the way toward the mean across agents for the same criterion. After three rounds, each score's distance from this mean is approximately 27\% of its initial distance. Invalid LLM responses receive up to two retries.
For the reliability measures, we report 95\% percentile bootstrap intervals based on $5{,}000$ resamples of the ten base descriptions.

To examine the scheduling consequences of the generated preferences,
we form the sixteen weight samples $\Lambda^{S,\pi}$ described in
Subsection~\ref{sec:explicit-uncertainty}. Each pre-deliberation and
mechanical sample contains thirty weight vectors, one for each agent
and base description. Each numerical and argumentative sample
contains $150$ vectors because it also includes the five deliberation
trajectories. For each sample, the mean vector and the uncertainty set
defined in~\eqref{eq:explicit-uncertainty-set} provide the inputs to
problems (MW) and (RB) in~\eqref{eq:explicit-problems}, respectively.

We solve both problems for the same synthetic clinic session, in
which one provider schedules $n=7$ patient requests over $H=300$
minutes, using $m=10$ appointment slots of length $\tau=30$ minutes.
The service durations and earliest acceptable slots are
\[
 d=(47,32,49,37,37,43,46),
 \qquad
 a=(0,0,0,1,1,2,3).
\]
For the objectives in~\eqref{eq:explicit-scheduling-objectives}, we set
$\nu_w=120$ and $\nu_u=\nu_o=60$ minutes. The larger reference scale
for total patient waiting time accounts for its accumulation across
patients and helps give the three normalized criteria comparable
magnitudes. These scales are part of the mapping from importance scores to objective weights, so the selected schedules depend on them as well as on the elicited scores. All problems are solved with HiGHS
\citep{huangfu_parallelizing_2018} through SciPy
\citep{virtanen_scipy_2020}.

%%%%%%%%%%%%%%%%%%%%%%%%%%%%%%%%%%%%%%%%%%%%%%%%%%%%%%%%%%
\subsubsection{Contextual sensitivity (Q1)}
\label{sec:explicit-context-results}

\begin{figure}[t]
 \centering
 \includegraphics[width=\textwidth]{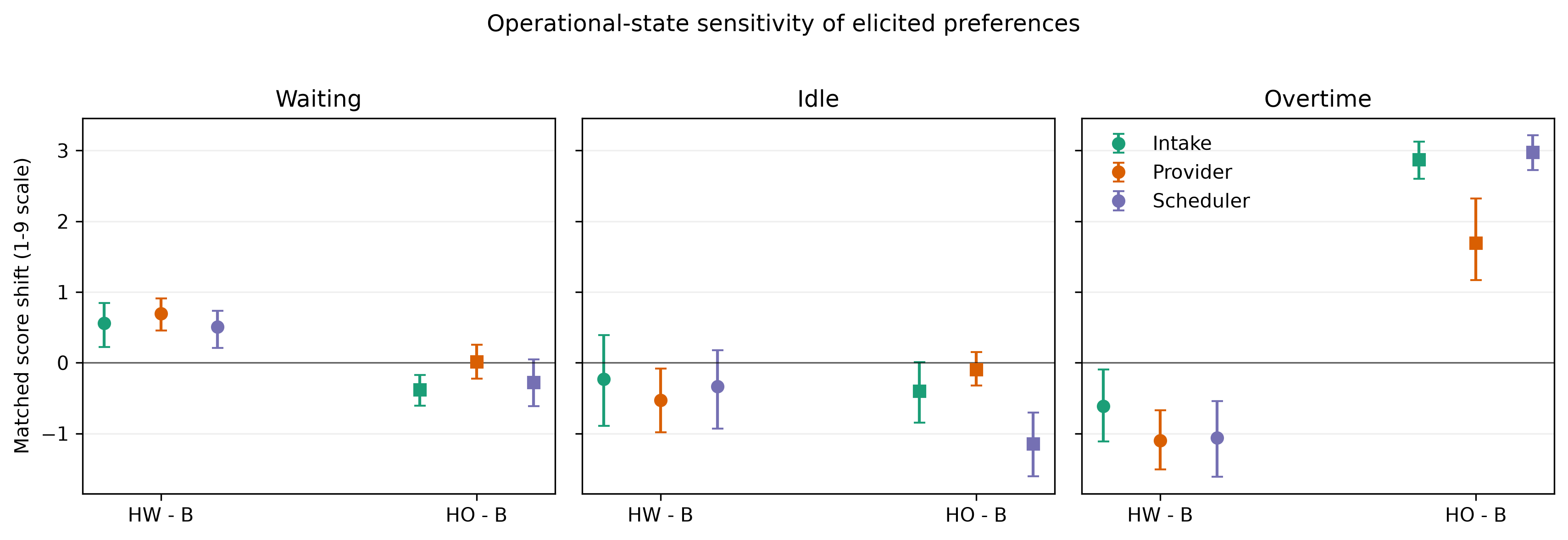}
 \caption{Paired score differences between condition HW or HO and condition B before deliberation, by agent and criterion. Values are mean differences in importance score points on the $1$--$9$ scale; positive values mean higher importance in the stated condition. Error bars are 95\% percentile bootstrap intervals over the ten base clinic descriptions.}
 \label{fig:explicit-contextual-shifts}
\end{figure}

Before deliberation, the paired differences~\eqref{eq:explicit-context-shift} show that the agents respond to the manipulated clinic conditions (Figure~\ref{fig:explicit-contextual-shifts}). Under HW, the mean paired difference in the waiting score is 0.560 points for Intake (95\% bootstrap interval $[0.225,0.845]$), 0.695 for Provider ($[0.455,0.910]$), and 0.506 for Scheduler ($[0.206,0.736]$). All three intervals lie above zero, supporting an increase in the importance assigned to waiting when the displayed mean waiting time rises. Under HO, the overtime score increases by 2.872, 1.691, and 2.977 points for Intake, Provider, and Scheduler, respectively; all three intervals exclude zero. In contrast, HW decreases the overtime score by 0.611, 1.097, and 1.062 points for the three agents.

\begin{figure}[t]
 \centering
 \includegraphics[width=0.90\textwidth]{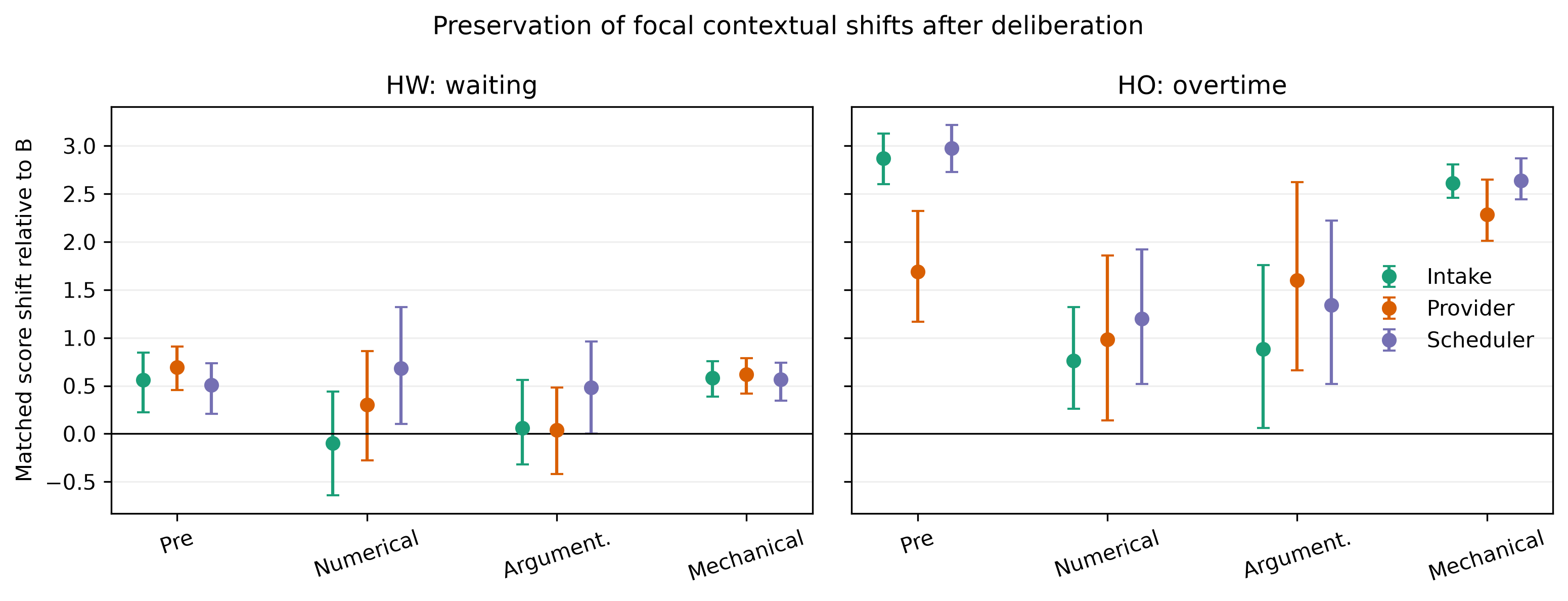}
 \caption{Focal paired score differences before deliberation and after round 3 of each protocol: waiting score under HW (left) and overtime score under HO (right), in importance score points on the $1$--$9$ scale. Error bars are 95\% percentile bootstrap intervals over the ten base clinic descriptions.}
 \label{fig:explicit-shift-preservation}
\end{figure}

Deliberation attenuates these responses by an amount that depends on the agent and criterion (Figure~\ref{fig:explicit-shift-preservation}). Averaged across agents, the HW--B difference in the waiting score decreases from 0.587 points before deliberation to 0.293 under the numerical protocol and 0.193 under the argumentative protocol. After numerical deliberation, the agent-specific mean differences are $-0.100$, 0.300, and 0.680 for Intake, Provider, and Scheduler, respectively; after argumentative deliberation, they are 0.060, 0.040, and 0.480. For the HO--B difference in the overtime score, the mean across agents decreases from 2.513 before deliberation to 0.980 under the numerical protocol and 1.273 under the argumentative protocol, although all agent-specific differences remain positive. The argumentative protocol preserves a larger mean overtime response than the numerical protocol but a smaller mean waiting response, so neither LLM protocol consistently preserves the initial contextual responses better than the other (\textbf{Q1 takeaway}).

The argumentative protocol also records the arguments and clinic-state fields used during deliberation. Among the 1{,}170 valid argumentative turns in rounds 2 and 3, 1{,}118 (95.56\%) cite a valid preceding argument from another agent (Table~\ref{tab:explicit-auditability}). Explicit references to the manipulated clinic indicators are less frequent. Mean waiting time is mentioned in 9.38\% of valid HW turns and 3.36\% of valid B turns, while projected overtime is mentioned in 13.53\% of valid HO turns and 10.54\% of valid B turns. Current time, which is also manipulated under HO, is referenced in 81.88\% of valid HO turns and 76.23\% of valid B turns. Valid terminal outputs account for 99.06\% of numerical turns and 98.06\% of argumentative turns.

\begin{table}[t]
\centering
\caption{Technical failures, invalid outputs, and references during deliberation. Entries report counts, with percentages in parentheses. Technical failures are unsuccessful requests. Invalid outputs are completed responses that fail JSON or field validation. Reference percentages are calculated using only valid argumentative turns.}
\label{tab:explicit-auditability}
\small
\resizebox{\textwidth}{!}{%
\begin{tabular}{llrr}
\hline
Measure & Scope & Numerical & Argumentative \\
\hline
Valid terminal outputs & All turns & 1783/1800 (99.06\%) & 1765/1800 (98.06\%) \\
Technical terminal failures & All turns & 17/1800 (0.94\%) & 24/1800 (1.33\%) \\
Schema-invalid terminal outputs & All turns & 0/1800 (0.00\%) & 11/1800 (0.61\%) \\
Technical failures across attempts & All attempts & 131/1923 (6.81\%) & 164/2005 (8.18\%) \\
Schema-invalid outputs across attempts & All attempts & 0/1923 (0.00\%) & 67/2005 (3.34\%) \\
Valid other-agent argument reference & Rounds 2--3 & --- & 1118/1170 (95.56\%) \\
Mean-waiting reference & HW & --- & 42/448 (9.38\%) \\
Mean-waiting reference & B & --- & 15/446 (3.36\%) \\
Projected-overtime reference & HO & --- & 59/436 (13.53\%) \\
Projected-overtime reference & B & --- & 47/446 (10.54\%) \\
Current-time reference & HO & --- & 357/436 (81.88\%) \\
Current-time reference & B & --- & 340/446 (76.23\%) \\
\hline
\end{tabular}}
\end{table}

%%%%%%%%%%%%%%%%%%%%%%%%%%%%%%%%%%%%%%%%%%%%%%%%%%%%%%%%%%
\subsubsection{Role differentiation and inter-agent dispersion (Q2 and Q3)}
\label{sec:explicit-role-results}

The initial scores also differentiate the three roles. All 12 pre-deliberation role-differentiation scores $\mathrm{RD}_{e}^{(i,S)}$ have 95\% intervals entirely above zero (Table~\ref{tab:role-persistence-summary}). This differentiation weakens substantially after the LLM exchanges. After numerical deliberation, three of the 12 intervals remain entirely above zero, three contain zero, and six lie entirely below zero. After argumentative deliberation, the corresponding counts are one, eight, and three. Mechanical averaging, in contrast, retains positive intervals for all 12 scores; this holds by construction, since the update in~\eqref{eq:mechanical-update} brings the agents' raw scores closer together while preserving their ordering. The LLM protocols thus weaken or reverse the initial role differences although the role descriptions are unchanged in every round (\textbf{Q2 takeaway}).

\begin{table}[t]
\centering
\caption{Role-differentiation scores before and after deliberation. Entries are counts of the 12 scores per protocol (three agents under four conditions) whose 95\% bootstrap interval lies entirely above zero (positive), contains zero, or lies entirely below zero (reversed). Intervals are not adjusted for multiple comparisons.}
\label{tab:role-persistence-summary}
\begin{tabular}{lrrr}
\hline
Protocol & Positive & Interval contains zero & Reversed \\
\hline
Pre-deliberation & 12 & 0 & 0 \\
Numerical & 3 & 3 & 6 \\
Argumentative & 1 & 8 & 3 \\
Mechanical & 12 & 0 & 0 \\
\hline
\end{tabular}
\end{table}

\label{sec:explicit-dispersion-results}

\begin{figure}[t]
 \centering
 \includegraphics[width=0.90\textwidth]{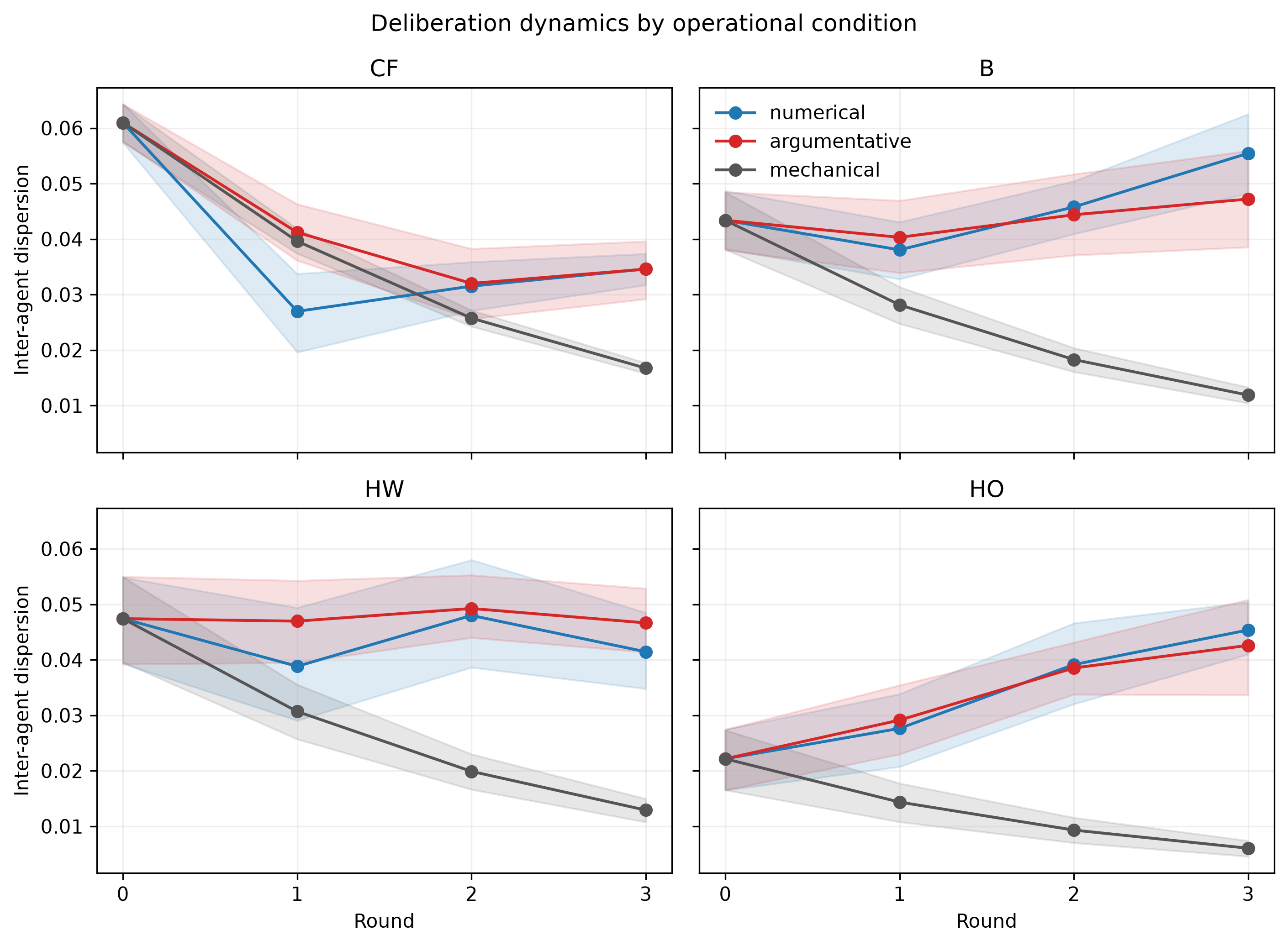}
 \caption{Inter-agent dispersion~\eqref{eq:explicit-dispersion} of the normalized weight vectors over the deliberation rounds, in units of normalized weight. Lines are means by protocol and condition; shaded regions are 95\% percentile bootstrap intervals over the ten base clinic descriptions.}
 \label{fig:explicit-dispersion}
\end{figure}

The LLM protocols do not reduce the inter-agent dispersion~\eqref{eq:explicit-dispersion} uniformly across conditions (Figure~\ref{fig:explicit-dispersion}). Under CF, the mean value of $D^{(r)}$ decreases from 0.0610 before deliberation to approximately 0.0346 under both LLM protocols. Under HW, it decreases from 0.0474 to 0.0414 under the numerical protocol and to 0.0467 under the argumentative protocol. Under HO, however, it increases from 0.0222 to 0.0454 and 0.0426, respectively. Mechanical averaging reduces the observed dispersion of the normalized weights to 27.2--27.4\% of its initial value. 
% For differences between unnormalized scores, the exact contraction factor is $(1-0.35)^3=0.274625$. 
% Under HO, the increase in dispersion occurs alongside weaker role differentiation, so greater disagreement does not imply stronger preservation of stakeholder priorities.
Therefore, when contextual information is provided, we can conclude that disagreement among agents remains high even after LLM deliberation (\textbf{Q3 takeaway}).

%%%%%%%%%%%%%%%%%%%%%%%%%%%%%%%%%%%%%%%%%%%%%%%%%%%%%%%%%%
\subsubsection{Weight uncertainty and scheduling consequences (Q4)}
\label{sec:explicit-decision-results}

The fixed service workload is $\sum_{p\in\mathcal P}d_p=291$ minutes, leaving $H-\sum_{p\in\mathcal P}d_p=9$ minutes of slack within the session horizon. We enumerate all 120 feasible schedules in~$X$ and calculate, for each schedule, its total patient waiting time, provider idle time, and overtime. These three totals form an outcome profile. Using the terminology in Section~\ref{subsec:background}, a profile is dominated if another feasible schedule gives equal or smaller values for all three totals and a strictly smaller value for at least one. Discarding dominated profiles and removing duplicates leaves five distinct nondominated profiles, expressed in minutes as
\[
 (\sum_{p\in\mathcal P}q_p,\sum_{p\in\mathcal P}g_p,O)
 \; \in \; \{
 (22,38,29),\,
 (23,25,16),\,
 (50,15,6),\,
 (72,8,0),\,
 (86,0,0)\}.
\]
For every $x\in X$, the completion time satisfies
$b_n+d_n=\sum_{p\in\mathcal P}d_p+\sum_{p\in\mathcal P}g_p$.
Consequently, the definition of overtime gives
$O=\max\{0,\sum_{p\in\mathcal P}g_p-9\}$.
Thus, up to nine minutes of provider idle time can be accommodated without overtime; any idle time beyond nine minutes produces overtime.
Every nondominated schedule has $x_1=0$ and hence no initial gap, so the idle time counted by~$f_u$ in~(MW) and (RB) consists entirely of gaps between patients.

We now examine the mean weights $\bar\lambda^{S,\pi}$ and uncertainty sets $\mathcal U_{S,\pi}$ obtained from the agents' scores after the final round of numerical or argumentative deliberation. The argumentative protocol produces wider componentwise intervals $[L_c^{S,\pi},U_c^{S,\pi}]$ than the numerical protocol for ten of the twelve criterion--condition combinations; the exceptions are the waiting-weight intervals under HW and HO (Table~\ref{tab:appendix-uncertainty} in Appendix~\ref{app:explicit-uncertainty}). Under the argumentative protocol ($\pi=\mathrm{arg}$), the mean weight vectors are
$\bar\lambda^{\mathrm{HW},\mathrm{arg}}=(0.493,0.184,0.323)$
and
$\bar\lambda^{\mathrm{HO},\mathrm{arg}}=(0.405,0.178,0.418)$.
Thus, the HO weights assign substantially more importance to overtime when evaluating the same feasible-schedule set $X$.

\begin{figure}[t]
 \centering
 \includegraphics[width=\textwidth]{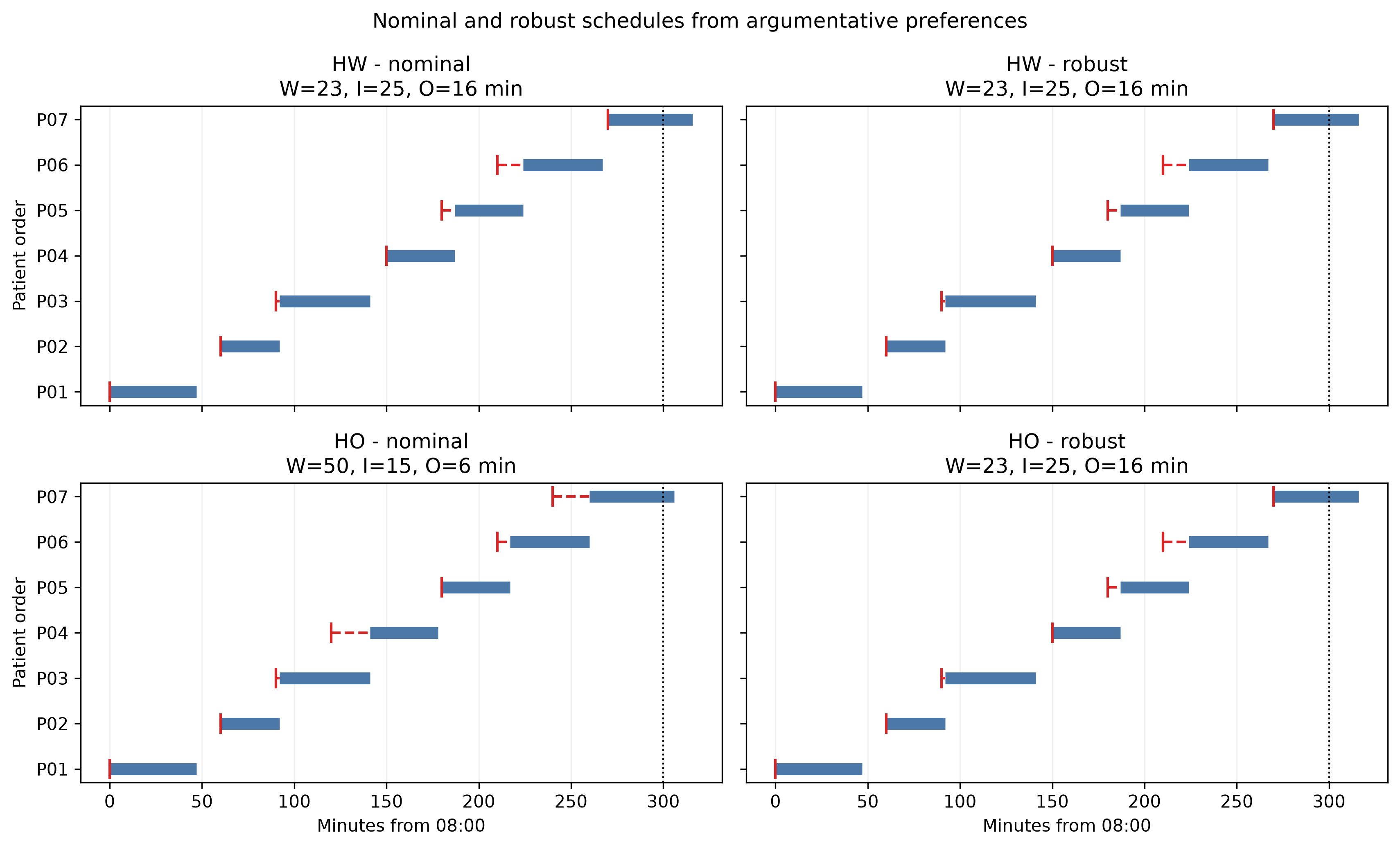}
 \caption{Schedules selected by problems (MW) and (RB) on the same instance with seven requests, using weights elicited under HW and HO and revised through argumentative deliberation. HW and HO label the elicitation conditions; all panels use the same session, requests, and constraints. Blue segments are service intervals, red markers are appointment times, and dashed red segments are patient waiting. Panel titles report total waiting (W), idle time (I), and overtime (O) in minutes.}
 \label{fig:explicit-schedules}
\end{figure}

Figure~\ref{fig:explicit-schedules} illustrates the resulting scheduling decisions. For $(S,\pi)=(\mathrm{HW},\mathrm{arg})$, both (MW) and (RB) select the nondominated profile $(23,25,16)$, corresponding to 23 minutes of total waiting, 25 minutes of idle time, and 16 minutes of overtime. For $(S,\pi)=(\mathrm{HO},\mathrm{arg})$, problem (MW) instead selects $(50,15,6)$, accepting 27 additional minutes of waiting in exchange for 10 fewer minutes of both idle time and overtime. Problem (RB), however, selects $(23,25,16)$. Evaluated at the mean weight vector, the robust solution increases the objective from 0.2548 to 0.2629 (3.18\%) relative to the (MW) solution, while reducing the worst-case objective over $\mathcal U_{\mathrm{HO},\mathrm{arg}}$ from 0.2959 to 0.2830 (4.37\%). Although overtime has the largest mean weight under HO, (RB) selects a schedule with more overtime than (MW). Across the 16 combinations of condition and weight source, (MW) and (RB) select different schedules in ten cases, with the resulting solutions spanning four distinct nondominated profiles (Appendix~\ref{app:explicit-schedules}). Accounting for uncertainty in the weights can therefore change the selected appointment schedule (\textbf{Q4 takeaway}).

%%%%%%%%%%%%%%%%%%%%%%%%%%%%%%%%%%%%%%%%%%%%%%%%%%%%%%%%%%
% SECTION 4
%%%%%%%%%%%%%%%%%%%%%%%%%%%%%%%%%%%%%%%%%%%%%%%%%%%%%%%%%%
\section{Implicit preference inference in emergency department planning}
\label{sec:implicit_preference}

The importance scores elicited in Section~\ref{sec:explicit-preferences} do not establish how an~LLM agent trades an improvement in one objective against a deterioration in another. We apply the LLM-PADA framework to infer this trade-off from an agent's choices between pairs of emergency department (ED) configurations. For each configuration, the agent sees two numerical outcomes: a measure of patient waiting time and the number of hours per week that a minor injuries unit (MIU) is open, where the MIU is a dedicated pathway for eligible low-acuity patients. Opening the~MIU for more hours can reduce waiting times but requires a greater commitment of ED resources. The choice of the agent may depend on these outcomes and on which configuration is presented as Option~1 (we refer to this as a \emph{presentation effect}). We estimate a trade-off weight using a choice model that includes both the outcomes and presentation order, and assess its predictive value against a benchmark that uses presentation order alone.

We ask LLM agents prompted from patient, physician, and operations perspectives to compare configurations reported in the simulation study of \citet{boresta_managing_2024}.

\begin{remark}
As in Section~\ref{sec:explicit-preferences}, we do not claim that the inferred weights represent the preferences of human patients, physicians, or hospital operators.
\end{remark}

%%%%%%%%%%%%%%%%%%%%%%%%%%%%%%%%%%%%%%%%%%%%%%%%%%%%%%%%%%
\subsection{Problem setting}
\label{sec:choice-setting}

We use the discrete-event simulation study of \citet{boresta_managing_2024} to define the operational setting. Its model represents stochastic patient arrivals, triage, and patient flow through the ED's Medical Unit (MU), Surgical Unit (SU), and MIU. Let $D=\{1,\ldots,7\}$ index the days of a week. A configuration
$x=(o_1,q_1,\ldots,o_7,q_7,\pi_{\mathrm m},\pi_{\mathrm a})$
specifies the MIU opening time $o_d$ and closing time $q_d$ on each day $d\in D$, together with the percentages $\pi_{\mathrm m}$ and $\pi_{\mathrm a}$ of eligible patients routed to the MIU in the morning and afternoon, respectively. If the MIU is closed on a day, it follows that~$o_d=q_d$. The number of MIU rooms is fixed in this case study.

\begin{figure}[tb]
 \centering
 \includegraphics[width=0.82\textwidth]{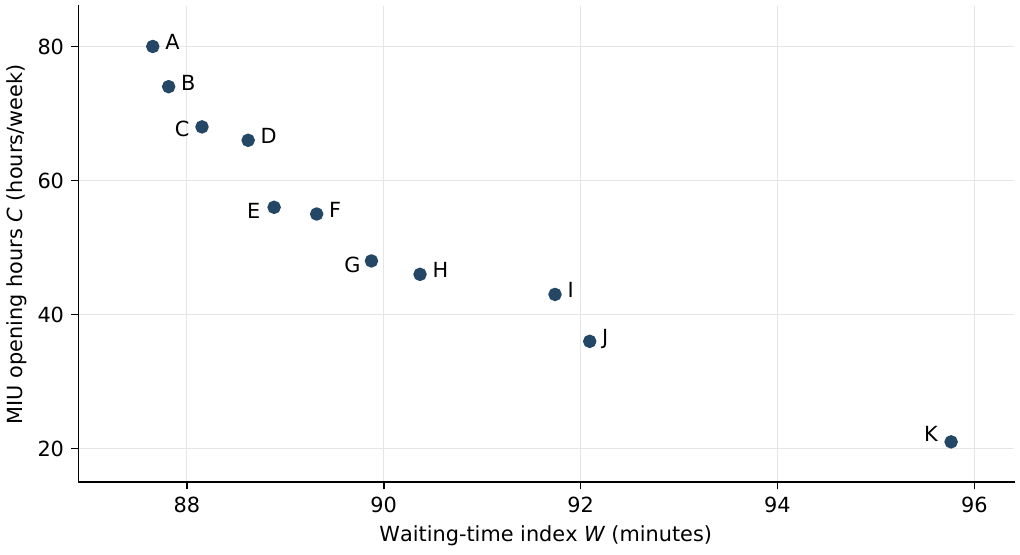}
 \caption{Average of the eight group mean waiting times, $f_w$ (minutes), and MIU opening hours, $f_{\mathrm{cost}}$ (hours per week), for the eleven configurations in Table~3 of \citet{boresta_managing_2024}.}
 \label{fig:choice-miu-tradeoffs}
\end{figure}

To define the average waiting time in the~ED, let $\mathcal T=\{R,Y,G,W\}$ denote the red, yellow, green, and white triage tags, and let~$U(\cdot)$ be a mapping that assigns to each triage tag the corresponding set of units. Specifically,
\[
 U(R)=U(Y)=\{\mathrm{MU},\mathrm{SU}\},\qquad
 U(G)=\{\mathrm{MU},\mathrm{SU},\mathrm{MIU}\},\qquad
 U(W)=\{\mathrm{MIU}\}.
\]
The set $\mathcal G=\{(t,u):t\in\mathcal T,\ u\in U(t)\}$ therefore contains eight patient groups. For $g\in\mathcal G$, let $T_g(x,\xi)$ be the simulated time from the start of triage to the first medical visit under configuration $x$, where $\xi$ represents the random inputs to the simulation. Define
\begin{equation}
 f_w(x)=\frac{1}{8}\sum_{g\in\mathcal G}
       \mathbb E_\xi[T_g(x,\xi)],
 \qquad
 f_{\mathrm{cost}}(x)=\sum_{d\in D}(q_d-o_d).
 \label{eq:choice-raw-objectives}
\end{equation}
The first quantity is the average of the eight group mean waiting times, in minutes. The second is the total number of MIU opening hours per week, used here as a resource measure rather than a monetary cost. The waiting-time objective reported in Table~3 of \citet{boresta_managing_2024} is the \emph{sum} of the eight group mean waiting times. Accordingly, we divide each waiting-time value by eight.
Now, let $\mathcal A=\{\mathrm A,\ldots,\mathrm K\}$ index the eleven configurations reported in Table~3 of \citet{boresta_managing_2024}, and let $x^a$ denote the complete configuration associated with label $a\in\mathcal A$. Our finite decision set is $X=\{x^a:a\in\mathcal A\}$. Thus, $x^a$ is an operating configuration, and the pair of objective values associated with it is $(f_w(x^a),f_{\mathrm{cost}}(x^a))$; see Figure~\ref{fig:choice-miu-tradeoffs}.

The eleven reported configurations are mutually nondominated with respect to the two reported objectives: none has both a lower value of $f_w$ and a lower value of $f_{\mathrm{cost}}$ than another. Selecting among them therefore requires a trade-off. We normalize the objectives and define
\begin{equation}
 \begin{aligned}
 w(x)&=\frac{f_w(x)-85}{15},&
 c(x)&=\frac{f_{\mathrm{cost}}(x)-21}{63},\\
 J_\lambda(x)&=\lambda w(x)+(1-\lambda)c(x),&
 \lambda&\in[0,1].
 \end{aligned}
 \label{eq:choice-objective}
\end{equation}
The normalization uses the intervals 85--100 minutes for $f_w$ and 21--84 hours per week for $f_{\mathrm{cost}}$. Thus, a waiting time of 85 minutes corresponds to $w(x)=0$ and a waiting time of 100 minutes to $w(x)=1$, while 21 and 84 opening hours per week correspond to $c(x)=0$ and $c(x)=1$, respectively. The reported values of~$f_w$ and~$f_{\mathrm{cost}}$ for all eleven configurations lie within these intervals. The weight $\lambda$ specifies the relative importance of the normalized waiting-time objective: $\lambda=1$ places all weight on waiting time, $\lambda=0$ places all weight on opening hours, and an interior value combines them. For a specified $\lambda$, we select a configuration by solving
\begin{equation}\label{prob:config_selection}
\min_{x\in X}J_\lambda(x).
\end{equation}

% Note that 
% weighted sums are a standard scalarization of multiobjective problems \citep{ehrgott_multicriteria_2005,miettinen_nonlinear_2012}. 
Over the finite set $X$, a weighted sum can select only
configurations whose normalized objective values lie on the lower
convex hull of $\{(w(x),c(x)):x\in X\}$. Nondominance alone
does not ensure that a configuration minimizes $J_\lambda$ for some
$\lambda$. Therefore, although configurations~D, F, H, and~I in
Figure~\ref{fig:choice-miu-tradeoffs} are nondominated, they do not
minimize $J_\lambda$ for any fixed $\lambda\in[0,1]$.
This restriction applies to problem~\eqref{prob:config_selection} with a fixed weight. It does not apply to the robust problem
in Section~\ref{sec:choice-optimization}, which can select configuration~H.

%%%%%%%%%%%%%%%%%%%%%%%%%%%%%%%%%%%%%%%%%%%%%%%%%%%%%%%%%%
\subsection{Research questions}
\label{sec:choice-questions}

The patient, physician, and operations prompts each describe the two outcomes, average waiting times and MIU opening hours, from their respective perspectives, but none assigns a numerical weight to either outcome. Unlike Section~\ref{sec:explicit-preferences}, where each role description is associated with one criterion, we do not posit an expected weight for any role and estimate the weight from the agent's pairwise choices. Our analysis asks:

\begin{itemize}
 \item[(Q1)] \emph{Predictive relevance:} For pairs of configurations not used to estimate the trade-off, do their waiting-time and opening-hour outcomes improve prediction of the agent's choice beyond a model that accounts for which option is presented first?
 \item[(Q2)] \emph{Choice consistency:} How consistent are the agents' choices across repeated requests and reversed presentation orders, and to what extent do they favor a configuration that is strictly better on both outcomes?
 \item[(Q3)] \emph{Weight uncertainty:}  How much uncertainty remains in the inferred trade-off weights across the tested models and stakeholder roles?
 \item[(Q4)] \emph{Configuration consequences:} How does
accounting for uncertainty in the estimated trade-off weight
affect the selected configuration compared with using a single weight?
\end{itemize}

%%%%%%%%%%%%%%%%%%%%%%%%%%%%%%%%%%%%%%%%%%%%%%%%%%%%%%%%%%
\subsection{Methodology}
\label{sec:choice-methodology}

We now describe how the four steps of LLM-PADA are implemented for ED configuration selection. Step~1 infers a trade-off weight from pairwise choices. Step~2 assesses predictive relevance and choice consistency to answer (Q1) and (Q2). Step~3 constructs an uncertainty set for the inferred weight and supports the analysis of (Q3). Step~4 compares mean-weight, fitted-weight, and robust (McRow) configuration selection to answer (Q4).

%%%%%%%%%%%%%%%%%%%%%%%%%%%%%%%%%%%%%%%%%%%%%%%%%%%%%%%%%%
\subsubsection{Preference generation and inference (Step 1)}
\label{sec:choice-design}

We first give each agent the definitions and units of $f_w$ and $f_{\mathrm{cost}}$ and the scaling intervals used for normalization in~\eqref{eq:choice-objective}. In a separate request for each comparison, we present two configurations in a two-row table containing their numerical values of $f_w$ and $f_{\mathrm{cost}}$. The rows are labeled Option~1 and Option~2. We present each pair in both orders and we ask each agent to select their preferred configuration. Below, we describe the details of the approach, and we use the term \emph{presentation effect} for any tendency to select the configuration shown as Option~1.

Consider an unordered pair $h=\{i,j\}$ with distinct $i,j\in\mathcal A$. Label its members so that~$f_w(x^i)<f_w(x^j)$ and $f_{\mathrm{cost}}(x^i)>f_{\mathrm{cost}}(x^j)$. Thus, $x^i$ denotes the configuration with lower waiting time but more MIU opening hours. Define $\Delta w_h=w(x^j)-w(x^i)>0$ and~$\Delta c_h=c(x^i)-c(x^j)>0$. Let $z=1$ when $x^i$ is presented as Option~1 and $z=-1$ when it is presented as Option~2. For each presentation order $z$, let $Y_{h,z}=1$ if the agent selects the lower-waiting configuration $x^i$, and $Y_{h,z}=0$ if it selects $x^j$.
For example, suppose $x^i$ has outcomes $(90\text{ minutes},60\text{ hours})$ and $x^j$ has outcomes $(94\text{ minutes},48\text{ hours})$. In the comparison $(h,1)$, $x^i$ is Option~1; selecting it gives $Y_{h,1}=1$. In the separate comparison $(h,-1)$, $x^j$ is Option~1; selecting it gives $Y_{h,-1}=0$.
The responses $Y_{h,1}$ and $Y_{h,-1}$ therefore correspond to separate comparisons of the same pair in opposite orders. For each combination of LLM model and agent role, we fit an \emph{outcome-based choice model}, which is a logistic regression that predicts selection of the lower-waiting configuration from the two normalized outcome differences and the presentation order:
\begin{equation}
 \begin{aligned}
 p_{h,z}
 &:=\Pr(Y_{h,z}=1\mid h,z)
   =\sigma\!\left(\theta_w\Delta w_h-\theta_c\Delta c_h+\gamma z\right),\quad\theta_w,\theta_c\geq 0,\quad
 \sigma(t)=\frac{1}{1+e^{-t}}.
 \end{aligned}
 \label{eq:choice-logit}
\end{equation}
The argument of $\sigma$ is the log-odds of choosing the lower-waiting configuration $x^i$. Increasing its waiting-time advantage $\Delta w_h$ raises those log-odds by $\theta_w$ per normalized unit; increasing its additional opening hours $\Delta c_h$ lowers them by $\theta_c$ per normalized unit. The sign constraints encode the direction in which each objective affects the fitted choice probability. The coefficient~$\gamma$ measures the presentation effect: when $\gamma>0$, presenting $x^i$ as Option~1 rather than Option~2 raises the log-odds of selecting it by $2\gamma$. When $\gamma=0$, $p_{h,1}=p_{h,-1}$ for every pair $h$.

Equation~\eqref{eq:choice-logit} models the probability of selecting $x^i$ using its waiting-time advantage $\Delta w_h$, its additional opening hours $\Delta c_h$, and its presentation order $z$. The same coefficients $\theta_w$, $\theta_c$, and $\gamma$ apply to every pair for a given agent role. It is related to the well-known Bradley--Terry model \citep{bradley_rank_1952}, with a term for the presentation effects documented in LLM comparisons \citep{ zheng_judging_2023}.
We now relate the choice model in~\eqref{eq:choice-logit} to the weighted objective in~\eqref{eq:choice-objective}. If $\kappa:=\theta_w+\theta_c>0$, define $\lambda=\theta_w/\kappa$. Then $\lambda\in[0,1]$, $\theta_w=\kappa\lambda$, and $\theta_c=\kappa(1-\lambda)$. Substitution gives
\[
 \theta_w\Delta w_h-\theta_c\Delta c_h
 =\kappa\bigl[\lambda\Delta w_h-(1-\lambda)\Delta c_h\bigr]
 =\kappa\bigl[J_\lambda(x^j)-J_\lambda(x^i)\bigr].
\]
The last equality follows from the definitions of $\Delta w_h$ and $\Delta c_h$. Since a smaller value of $J_\lambda$ is preferred, a positive difference $J_\lambda(x^j)-J_\lambda(x^i)$ increases the fitted probability of choosing $x^i$. The ratio $\lambda$ specifies the trade-off between waiting time and opening hours, while $\kappa$ determines how strongly differences in $J_\lambda$ affect choice probabilities.

Let $\mathcal P=\{\{i,j\}:i,j\in\mathcal A,\ i\ne j\}$ denote all unordered pairs of distinct configurations, and let $\mathcal H_{\mathrm{est}}\subset\mathcal P$ be the pairs used for estimation. For each pair, we request a choice in both orders. The estimation observations are indexed by $\mathcal D_{\mathrm{est}}=\{(h,z):h\in\mathcal H_{\mathrm{est}},\ z\in\{-1,1\}\}$; $y_{h,z}\in\{0,1\}$ denotes the observed choice. We describe how $\mathcal H_{\mathrm{est}}$ is selected below. We estimate $\theta_w$, $\theta_c$, and $\gamma$ by solving
\begin{equation}
 \begin{aligned}
 \min_{\theta_w,\theta_c,\gamma}\quad&
 \sum_{(h,z)\in\mathcal D_{\mathrm{est}}}
 \left[
   -y_{h,z}\ln p_{h,z}
   -(1-y_{h,z})\ln(1-p_{h,z})
 \right]
 +\frac{\rho}{2}
  \left(\theta_w^2+\theta_c^2+\gamma^2\right)\\
 \mathrm{s.t.}\quad&
 0\leq\theta_w,\theta_c\leq100,\qquad
 -10\leq\gamma\leq10,
 \end{aligned}
 \label{eq:choice-estimation}
\end{equation}
where $p_{h,z}$ is given by~\eqref{eq:choice-logit} and $\rho=0.01$. The first term is the negative Bernoulli log-likelihood, treating responses to separate requests as conditionally independent. The quadratic penalty discourages extreme coefficient values.
After solving~\eqref{eq:choice-estimation}, we compute $\widehat\lambda=\widehat\theta_w/(\widehat\theta_w+\widehat\theta_c)$ when the denominator exceeds $10^{-6}$, where~$\widehat\theta_w$ and~$\widehat\theta_c$ are the optimal solutions of~\eqref{eq:choice-estimation}. If both objective coefficients are effectively zero, their ratio does not define a meaningful trade-off weight. This does not occur in the reported fits.
Before collecting responses, we select a subset
$\mathcal H_{\mathrm{est}}\subset\mathcal P$ for estimation
and reserve the remaining pairs,
$\mathcal H=\mathcal P\setminus\mathcal H_{\mathrm{est}}$,
for holdout evaluation.

%%%%%%%%%%%%%%%%%%%%%%%%%%%%%%%%%%%%%%%%%%%%%%%%%%%%%%%%%%
\subsubsection{Reliability assessment (Step 2)}
\label{sec:choice-reliability}

We compare~\eqref{eq:choice-logit} with two benchmarks that use different information to predict a choice. The \emph{position-only benchmark}, $p^{\mathrm{pos}}_{h,z}=\sigma(\gamma_0z)$, assigns the same probability $\sigma(\gamma_0)$ of selecting Option~1 in every comparison, regardless of the configurations' waiting times and MIU opening hours. It can capture a tendency to choose Option~1, but does not use the numerical differences in waiting time or MIU opening hours.
The \emph{Bradley--Terry benchmark} \citep{bradley_rank_1952} is~$p^{\mathrm{BT}}_{h,z}
   =\sigma(v_i-v_j+\gamma_{\mathrm{BT}}z)$,
where $v_a$ is a score estimated for configuration $a\in\mathcal A$. A larger score makes a configuration more likely to be selected. This benchmark also accounts for presentation order, but it estimates a separate score for each configuration instead of requiring its choices to follow one common weighted trade-off between waiting time and opening hours. We compare this benchmark with~\eqref{eq:choice-logit} on pairs that were not used to fit either model, to assess whether estimating a separate score for each configuration improves prediction.

We fit the position-only benchmark by minimizing the same negative log-likelihood with $p_{h,z}$ replaced by $p^{\mathrm{pos}}_{h,z}$, optimizing over $-10\leq\gamma_0\leq10$, and adding the penalty $\rho\gamma_0^2/2$. For the Bradley--Terry benchmark, we use $p^{\mathrm{BT}}_{h,z}$, optimize over $-10\leq v_a,\gamma_{\mathrm{BT}}\leq10$ for all $a\in\mathcal A$, and add the penalty
$\rho(\sum_{a\in\mathcal A}v_a^2+\gamma_{\mathrm{BT}}^2)/2$.

After fitting the outcome-based model and the two benchmarks,
we hold their coefficients fixed and ask the agent to make
a new choice for every pair in $\mathcal P$, in both presentation
orders. This gives 110 evaluation responses per LLM model
and role. The 62 responses for the 31 holdout pairs in
$\mathcal H$ are used to assess prediction on pairs excluded
from estimation. The remaining 48 responses repeat the 24
estimation pairs in the same presentation orders as before.
Comparing these new choices with the original choices measures
repeatability when presentation order is unchanged. Let $m\in\{\mathrm{pos},\mathrm{out},\mathrm{BT}\}$ identify the position-only benchmark, outcome-based model, and Bradley--Terry benchmark, respectively. For each holdout comparison $(h,z)$, let $\widehat p^m_{h,z}$ denote model $m$'s fitted probability of selecting $x^i$ and let $y^{\mathrm{eval}}_{h,z}\in\{0,1\}$ denote the new observed choice. We compute the well-known Brier score \citep{brier_verification_1950} and the improvement over the position-only benchmark:
\begin{equation}
 B_m=\frac{1}{2|\mathcal H|}
 \sum_{h\in\mathcal H}\sum_{z\in\{-1,1\}}
 \left(y^{\mathrm{eval}}_{h,z}-\widehat p^m_{h,z}\right)^2,
 \qquad
 \Delta B=B_{\mathrm{pos}}-B_{\mathrm{out}}.
 \label{eq:choice-prediction}
\end{equation}
The Brier score in~\eqref{eq:choice-prediction} averages the squared difference between the observed binary choice $y^{\mathrm{eval}}_{h,z}$ and its predicted probability $\widehat p^m_{h,z}$ over the 62 holdout responses; smaller values indicate better predictions. Thus, $\Delta B>0$ means that the outcome-based model has a lower Brier score than the position-only benchmark. 
Let $\mathbf{1}\{A\}$ denote the indicator of a condition $A$: it equals 1 when $A$ is true and 0 otherwise. 
At a threshold of 0.5, the binary prediction $\mathbf{1}\{\widehat p^m_{h,z}\geq 0.5\}$ equals 1 when model $m$ predicts selection of the lower-waiting configuration and 0 otherwise. Two models can have identical binary predictions on every holdout comparison yet different Brier scores because their predicted probabilities differ.

We also report classification accuracy, the proportion of holdout responses correctly predicted after converting each fitted probability to a binary prediction at threshold 0.5. For $m\in\{\mathrm{pos},\mathrm{out},\mathrm{BT}\}$, using the predicted probabilities~$\widehat p^m_{h,z}$ and observed choices~$y^{\mathrm{eval}}_{h,z}$ from~\eqref{eq:choice-prediction}, define
\begin{equation}
 \widehat y^m_{h,z}=\mathbf 1\!\left\{\widehat p^m_{h,z}\geq0.5\right\},
 \qquad
 A_m=\frac{1}{2|\mathcal H|}
 \sum_{h\in\mathcal H}\sum_{z\in\{-1,1\}}
 \mathbf 1\!\left\{\widehat y^m_{h,z}=y^{\mathrm{eval}}_{h,z}\right\}.
 \label{eq:choice-accuracy}
\end{equation}
Here, $\widehat y^m_{h,z}=1$ predicts selection of the lower-waiting configuration, and $A_m$ averages correct predictions over the 62 holdout responses. Accuracy expressed as a percentage is $100A_m$.

To assess how $\Delta B$ varies across the holdout pairs, we resample the 31 pairs with replacement 5{,}000 times, keeping the two presentation orders of each pair together. We recompute $\Delta B$ for each resample without refitting the models and report its 2.5th and 97.5th percentiles. This interval reflects variation across the selected holdout pairs.

We also report two behavioral diagnostics. For each $h\in\mathcal H$, let $a_h^{(z)}\in\mathcal A$ be the configuration selected in presentation order $z$. \emph{Order agreement} is the proportion of holdout pairs for which the agent selects the same configuration in both presentation orders:
\[
 \mathrm{OA}=\frac{1}{|\mathcal H|}
 \sum_{h\in\mathcal H}
 \mathbf 1\!\left[a_h^{(1)}=a_h^{(-1)}\right].
\]
An agent that always selects Option~1 has $\mathrm{OA}=0$, while an agent that makes the same deterministic choice regardless of order has $\mathrm{OA}=1$. Even if $\gamma=0$, independent stochastic responses can disagree across orders.

All comparisons among the eleven simulated ED configurations involve a trade-off: the configuration with lower waiting time has more MIU opening hours. To assess whether presentation order can lead an agent to choose an alternative that is worse on both outcomes, we construct five additional pairs specified by numerical waiting times and opening hours rather than drawn from the eleven configurations. In each pair, one alternative has both lower waiting time and fewer opening hours than the other, so it is preferred under $J_\lambda$ for every $\lambda\in[0,1]$. We present each pair in both orders, yielding ten additional choices per model and role. We define as \textit{dominance compliance} the number of these ten choices in which the agent selects the alternative with both lower outcomes. These responses are excluded from model estimation and holdout prediction.

We summarize evidence for a fitted trade-off through an exploratory \emph{diagnostic screen}. A combination of model and role meets the screen if the lower endpoint of its $\Delta B$ interval is positive, dominance compliance is 10/10, and $0<\widehat\lambda<1$. The first two conditions require evidence of predictive improvement beyond presentation order and observed compliance with strict dominance. The third identifies a fitted trade-off involving both objectives: since $\widehat\lambda=\widehat\theta_w/(\widehat\theta_w+\widehat\theta_c)$, an interior weight requires both objective coefficients to be positive. An endpoint weight can represent exclusive priority for one objective and is not intrinsically invalid; the screen is specifically intended to summarize evidence for a trade-off between the two objectives. We report order agreement separately. The screen was defined after inspecting the results and is descriptive; its intervals are not adjusted for multiple comparisons.

%%%%%%%%%%%%%%%%%%%%%%%%%%%%%%%%%%%%%%%%%%%%%%%%%%%%%%%%%%
\subsubsection{Weight uncertainty construction (Step 3)}
\label{sec:choice-decisions}
\label{sec:choice-uncertainty}

The bootstrap interval quantifies estimation uncertainty under the outcome-based model. To construct it, we use a parametric bootstrap~\citep{efron_bootstrap_1979}: we generate 1{,}000 synthetic estimation datasets by independently sampling a binary response from the fitted Bernoulli probability for each ordered estimation comparison. For each dataset, we solve~\eqref{eq:choice-estimation} again with the same bounds and penalty, jointly reestimating the waiting-time, opening-hour, and presentation coefficients. We retain a refit when the numerical solver reports successful convergence for this optimization problem and the sum of its two objective coefficients exceeds $10^{-6}$. This second condition ensures that the ratio defining the weight is numerically meaningful. Indexing the retained refits by $b=1,\ldots,K$, their weights are
\begin{equation}
 \widehat\lambda^{(b)}=
 \frac{\widehat\theta_w^{(b)}}
 {\widehat\theta_w^{(b)}+\widehat\theta_c^{(b)}}.
 \label{eq:choice-bootstrap-weight}
\end{equation}

The endpoints $L$ and $U$ are the empirical 2.5th and 97.5th percentiles of these $K$ retained weights. These endpoints define a 95\% bootstrap interval for each model and role. Thus, the interval is based on successful refits with defined weights, rather than necessarily on all 1{,}000 generated datasets.

For each combination of model and role, we use the $K$ retained bootstrap weights defined in~\eqref{eq:choice-bootstrap-weight} to construct
\[
 \bar\lambda=\frac{1}{K}\sum_{b=1}^{K}\widehat\lambda^{(b)},
 \qquad \mathcal U_\lambda=[L,U],
\]
where $L$ and $U$ are the percentile endpoints defined above. We adopt $\mathcal U_\lambda$ as the weight uncertainty set. The 95\% level describes the bootstrap interval construction; it does not provide a guarantee that the inferred weight represents a human stakeholder's preference or that the selected configuration improves operational performance.

This interval describes uncertainty in the estimated weight conditional on the fitted probability model and estimation procedure. It does not account for departures from the assumption that one common trade-off applies across pairs, nor does it establish that removing the fitted presentation term would recover an order-independent preference. We use it below to examine how uncertainty in this model-based estimate affects configuration selection.

To answer (Q3), we examine the fitted weights and their bootstrap intervals across the six combinations of model and role, and compare the fitted estimates with the bootstrap means.

%%%%%%%%%%%%%%%%%%%%%%%%%%%%%%%%%%%%%%%%%%%%%%%%%%%%%%%%%%
\subsubsection{Configuration decision analysis (Step 4)}
\label{sec:choice-optimization}

Step~4 uses the inferred weights and their uncertainty set to address (Q4). Throughout this analysis, the reported simulation estimates of $f_w(x)$ and the opening-hour values $f_{\mathrm{cost}}(x)$ are held fixed for every $x\in X$. The uncertainty set therefore concerns the inferred weight and does not incorporate uncertainty in the waiting-time estimates.

Using $J_\lambda$ from~\eqref{eq:choice-objective}, we compare three decision rules:
\begin{equation}
 \begin{aligned}
 \text{(MW)}\quad &x^{\mathrm{MW}}\in\operatorname*{arg\,min}_{x\in X}J_{\bar\lambda}(x),\\
 \text{(FW)}\quad &x^{\mathrm{FW}}\in\operatorname*{arg\,min}_{x\in X}J_{\widehat\lambda}(x),\\
 \text{(RB)}\quad &x^{\mathrm{RW}}\in\operatorname*{arg\,min}_{x\in X}R(x),
 \qquad R(x):=\max_{\lambda\in\mathcal U_\lambda}J_\lambda(x).
 \end{aligned}
 \label{eq:choice-problems}
\end{equation}
The mean-weight rule (MW) uses the bootstrap mean, whereas the fitted-weight rule (FW) uses the estimate obtained from the observed estimation responses. We retain both because the bootstrap mean can differ from the fitted weight, particularly near a coefficient bound. Since $J_\lambda(x)$ is affine in $\lambda$,
\[
 J_{\bar\lambda}(x)=\frac{1}{K}\sum_{b=1}^{K}J_{\widehat\lambda^{(b)}}(x)
 \qquad\text{for every }x\in X.
\]
Hence, (MW) minimizes the weighted objective averaged over the retained bootstrap estimates. This identity provides a computational interpretation of the averaging rule; the bootstrap samples are not observations of a population of stakeholder preferences.

The robust rule (RB) minimizes the largest weighted objective over $\mathcal U_\lambda$. For any fixed configuration,
\begin{equation}
 R(x)=\max\{J_L(x),J_U(x)\},
 \label{eq:choice-robust-endpoints}
\end{equation}
because $J_\lambda(x)=c(x)+\lambda[w(x)-c(x)]$ is affine in $\lambda$, and an affine function attains its maximum over a closed interval at an endpoint. Equivalently, (RB) has the epigraph formulation
\[
 \begin{aligned}
 \min_{x\in X,\ t\in\mathbb R}\quad &t\\
 \mathrm{s.t.}\quad &t\geq Lw(x)+(1-L)c(x),\\
 &t\geq Uw(x)+(1-U)c(x).
 \end{aligned}
\]
All three problems are solved exactly by enumerating the eleven configurations. Unlike (MW) and (FW), (RB) minimizes a maximum of weighted sums. Over a finite decision set, it can therefore select a nondominated configuration that does not minimize any single weighted sum; configuration~H below illustrates this distinction. Such configurations must remain in the candidate set for the robust problem.

These decisions depend on the normalization in~\eqref{eq:choice-objective}. At a fixed weight, subtracting a constant from either normalized objective adds a configuration-independent constant to $J_\lambda$ and leaves its minimizers unchanged. Across weights, however, this constant generally depends on $\lambda$, so the worst-case objective and its minimizer can change. We therefore keep the reference origins and scale widths fixed when constructing and comparing the three rules. The robust criterion is a worst-case objective value, rather than regret relative to the optimal configuration at each weight.

To express the operational consequences of robustness relative to (MW), define
\[
 \Delta W=f_w(x^{\mathrm{RW}})-f_w(x^{\mathrm{MW}}),
 \qquad
 \Delta C=f_{\mathrm{cost}}(x^{\mathrm{RW}})-f_{\mathrm{cost}}(x^{\mathrm{MW}}).
\]
Thus, $\Delta W$ is the change in the average of the eight group mean waiting times, in minutes, and $\Delta C$ is the change in weekly MIU opening hours. We also report the loss at the mean weight and the reduction in the worst-case objective:
\begin{equation}
 \begin{aligned}
 P&=100\frac{J_{\bar\lambda}(x^{\mathrm{RW}})-J_{\bar\lambda}(x^{\mathrm{MW}})}
 {J_{\bar\lambda}(x^{\mathrm{MW}})},\\
 G&=100\frac{R(x^{\mathrm{MW}})-R(x^{\mathrm{RW}})}{R(x^{\mathrm{MW}})}.
 \end{aligned}
 \label{eq:choice-robust-comparison}
\end{equation}
The denominators are positive in all reported comparisons. Optimality of $x^{\mathrm{MW}}$ and $x^{\mathrm{RW}}$ implies $P\geq0$ and $G\geq0$. The quantity $P$ measures the increase in the objective at the mean weight incurred by selecting the robust configuration; $G$ measures the corresponding reduction in the worst-case objective. Both percentages depend on the normalized objective values and should be interpreted alongside the changes in minutes and opening hours.

%%%%%%%%%%%%%%%%%%%%%%%%%%%%%%%%%%%%%%%%%%%%%%%%%%%%%%%%%%
\subsection{Numerical results}
\label{sec:choice-numerical-results}

The experiments use locally executed Qwen3.5-9B and Qwen3.5-27B models with Q4\_K\_M quantization, temperature $0.7$, top-$p$ $0.9$, top-$k$ $20$, and thinking mode disabled. The 27B experiment was an exploratory replication conducted after examination of the 9B diagnostics; the elicitation and estimation procedures were unchanged.

The eleven configurations in Figure~\ref{fig:choice-miu-tradeoffs} give $|\mathcal P|=\binom{11}{2}=55$ pairs. Before collecting responses, we choose 24 pairs for $\mathcal H_{\mathrm{est}}$ and reserve the remaining 31 pairs, $\mathcal H=\mathcal P\setminus\mathcal H_{\mathrm{est}}$, for holdout evaluation. First, we choose eleven pairs that connect the configurations in a closed chain: each configuration is paired with the next, and the last is paired with the first. This ensures that every configuration is connected to every other through a sequence of selected comparisons. We then add thirteen pairs to cover different magnitudes and balances of the two outcome differences. For each candidate pair $h$, the weight $\lambda_h^*=\Delta c_h/(\Delta w_h+\Delta c_h)$ makes $\lambda_h^*\Delta w_h=(1-\lambda_h^*)\Delta c_h$, and hence $J_{\lambda_h^*}(x^i)=J_{\lambda_h^*}(x^j)$. It is the break-even weight: a weight above it favors the lower-waiting configuration $x^i$, whereas a weight below it favors the configuration with fewer opening hours $x^j$. We represent each candidate by $\Delta w_h$, $\Delta c_h$, and $\lambda_h^*$; after scaling these three coordinates, we repeatedly add the candidate farthest from its closest already selected pair.

For each Qwen model, we collect responses from the three agent roles. Each role provides 48 estimation responses (24 pairs in two orders), 110 evaluation responses (all 55 configuration pairs in two orders), and 10 dominance-check responses (five additional pairs in two orders). Thus, each Qwen model contributes $3(48+110+10)=504$ analyzed responses.

%%%%%%%%%%%%%%%%%%%%%%%%%%%%%%%%%%%%%%%%%%%%%%%%%%%%%%%%%%
\subsubsection{Predictive performance and choice consistency (Q1--Q2)}
\label{sec:choice-results}

Table~\ref{tab:choice-main} reports the fitted waiting-time weights and the measures defined in Sections~\ref{sec:choice-reliability} and~\ref{sec:choice-uncertainty}. To address (Q1), we compare the outcome-based model in~\eqref{eq:choice-logit} with the position-only benchmark using the holdout Brier-score improvement $\Delta B$ in~\eqref{eq:choice-prediction}. To address (Q2), we examine order agreement, which records whether the same configuration is selected in both presentation orders, and dominance compliance, which counts selections of the alternative with both lower waiting time and fewer opening hours in the ten additional comparisons. We also examine repeatability under a fixed presentation.

% Numerical entries retained from sec_4_codex_18set/09_export_paper_section4.py.
\begin{table}[tb]
 \centering
 \caption{Estimated waiting-time weights and predictive and behavioral diagnostics. Brackets give 95\% bootstrap intervals.}
 \label{tab:choice-main}
 \small
 \setlength{\tabcolsep}{4pt}
 \resizebox{\textwidth}{!}{%
 \begin{tabular}{llccccc}
 \hline
 Model & Role & $\widehat{\lambda}$ & \shortstack{$\Delta B$\\vs.\ position} & \shortstack{Order\\agreement} & \shortstack{Dominance\\compliance} & \shortstack{Diagnostic\\screen} \\
 \hline
 9B & Patient & $0.70\;[0.00,\,1.00]$ & $+0.006\;[-0.003,\,0.018]$ & 0.26 & 7/10 & Not met \\
 9B & Physician & $1.00\;[0.00,\,1.00]$ & $+0.002\;[-0.001,\,0.006]$ & 0.19 & 9/10 & Not met \\
 9B & Operations & $1.00\;[0.55,\,1.00]$ & $+0.002\;[-0.011,\,0.017]$ & 0.23 & 8/10 & Not met \\
 \hline
 27B & Patient & $0.75\;[0.71,\,1.00]$ & $+0.075\;[0.032,\,0.119]$ & 0.35 & 10/10 & Met \\
 27B & Physician & $0.72\;[0.67,\,1.00]$ & $+0.009\;[-0.030,\,0.048]$ & 0.29 & 10/10 & Not met \\
 27B & Operations & $0.64\;[0.58,\,0.80]$ & $+0.071\;[0.033,\,0.109]$ & 0.61 & 10/10 & Met \\
 \hline
 \end{tabular}%
 }
\end{table}

For Qwen3.5-9B, the outcome-based model in~\eqref{eq:choice-logit} attains holdout accuracies of 87.1\%, 90.3\%, and 85.5\% for Patient, Physician, and Operations, respectively. The position-only benchmark defined in Section~\ref{sec:choice-reliability} attains the same accuracies: both models predict selection of Option~1 in every holdout comparison. These accuracies therefore reflect the frequency of first-option selections and do not show an advantage from using outcome differences. The Brier-score improvements $\Delta B$ defined in~\eqref{eq:choice-prediction} range from 0.002 to 0.006, and all three intervals contain zero (Table~\ref{tab:choice-main}). Thus, the holdout results provide limited evidence of predictive improvement beyond presentation order for this model.

For Qwen3.5-27B, the results differ by role. For Patient, the Brier score decreases from 0.170 under the position-only benchmark to 0.095 under the outcome-based model, giving $\Delta B=0.075$ with interval $[0.032,0.119]$. For Operations, it decreases from 0.213 to 0.142, giving $\Delta B=0.071$ with interval $[0.033,0.109]$. Both intervals have positive lower endpoints, supporting a predictive contribution from the outcome differences. For Physician, $\Delta B=0.009$ with interval $[-0.030,0.048]$, so the same conclusion is not supported. Taken together, the holdout results support a predictive contribution from outcome differences beyond presentation order for 27B Patient and Operations, while remaining inconclusive for the other four combinations of model and role (\textbf{Q1 takeaway}).
The Bradley--Terry benchmark defined in Section~\ref{sec:choice-reliability} does not consistently improve the holdout Brier score over the outcome-based model. Assigning a separate score to each configuration therefore provides no consistent predictive advantage in this experiment.

For (Q2), the 9B behavioral diagnostics reinforce this limitation. Order agreement ranges from 0.19 to 0.26, and the alternative with both lower outcomes is selected in seven to nine of the ten dominance comparisons. All six selections of the dominated alternative across the three roles occur when the dominating alternative is presented as Option~2.

All three 27B roles select the dominating alternative in all ten additional comparisons, but their choices among the eleven nondominated configurations remain sensitive to presentation order. Option~1 is selected in 82.3\%, 82.3\%, and 69.4\% of the holdout responses for Patient, Physician, and Operations, respectively. The paired presentation design provides a 50\% reference rate: if an agent deterministically selects the same configuration in both orders, that configuration is Option~1 in exactly one of the two comparisons. If choices are stochastic but their probabilities do not depend on order, the expected proportion of Option~1 selections is likewise 50\%, although a realized sample need not attain this value. The observed proportions are above this reference and describe a tendency to select the first option. 

Order agreement is 0.35, 0.29, and 0.61 for the three 27B roles. Because each order is evaluated through a separate stochastic request, disagreement can reflect sampling variability as well as presentation effects; order agreement alone does not separate them. Across all six combinations of model and role, agreement between the original and repeated responses to the same ordered estimation comparison ranges from 75.0\% to 89.6\%. Thus, neither repeatability under a fixed presentation nor full compliance in the dominance comparisons establishes consistency across presentation orders when the alternatives involve a trade-off (\textbf{Q2 takeaway}).

Only 27B Patient and 27B Operations meet the screen defined in Section~\ref{sec:choice-reliability}. For 27B Physician, dominance compliance of 10/10 and an interior fitted weight coexist with an inconclusive holdout Brier comparison. None of the 9B combinations meets the screen. 

One possible explanation for the weaker Physician result concerns how the role description relates to the operational criteria. Patient and Operations perspectives can suggest concerns about waiting time and resource use, respectively, whereas a Physician perspective may give less direct guidance about their relative importance. However, none of the prompts specifies an objective weight, and the current experiment does not test this explanation. In particular, an inconclusive Brier-score comparison does not establish that the Physician agent is indifferent between the objectives or has no preference. As shown in Section~\ref{sec:choice-design}, $\lambda$ describes the relative coefficients, whereas $\kappa=\theta_w+\theta_c$ determines how strongly objective differences affect choice probabilities. A model can have an interior weight without yielding a detectable predictive improvement beyond presentation order. Testing the proposed role interpretation would require a controlled comparison of prompts that differ in how explicitly they connect the Physician perspective to these criteria.

%%%%%%%%%%%%%%%%%%%%%%%%%%%%%%%%%%%%%%%%%%%%%%%%%%%%%%%%%%
\subsubsection{Uncertainty in the inferred weights (Q3)}
\label{sec:choice-weight-results}

Table~\ref{tab:choice-main} reports the fitted weights and their
bootstrap intervals. For Qwen3.5-9B, the intervals for Patient and
Physician span $[0,1]$, while the Physician and Operations fits have
$\widehat\theta_c=0$ and hence $\widehat\lambda=1$. A zero opening-hour
coefficient means that the fitted model assigns no marginal effect
to additional opening hours after accounting for the waiting-time
difference and presentation order. Together with the holdout and
behavioral diagnostics, the broad intervals and endpoint estimates
provide weak support for interpreting the fitted ratios as
operational trade-offs.

For Qwen3.5-27B, the fitted weights are 0.75, 0.72, and 0.64 for
Patient, Physician, and Operations, respectively. Under the fixed
normalization in~\eqref{eq:choice-objective}, each point estimate
places more weight on waiting time than on opening hours.
The bootstrap intervals are narrower than their 9B counterparts
for all three roles, although the Patient and Physician intervals
still extend to $1$. Only the Operations interval is contained
in $(0,1)$. The intervals overlap substantially, and we do not
perform a formal comparison of the weights across roles.
The point estimates therefore describe the fitted objectives
without establishing differences between stakeholder roles.
The intervals indicate less uncertainty in the inferred weights
for 27B than for 9B, but only 27B Operations has an interval that
excludes both endpoint weights (\textbf{Q3 takeaway}).

The bootstrap intervals and means use the retained refits with
defined weights: $K$ is 991, 961, and 992 for 9B Patient, Physician,
and Operations, respectively, and 1{,}000 for each 27B role.
The resulting bootstrap means in Table~\ref{tab:choice-decisions}
can differ from the fitted estimates in Table~\ref{tab:choice-main}.
For 9B Physician and Operations, the fitted weights are $1$,
whereas the bootstrap means are $0.719$ and $0.896$, respectively.
Section~\ref{sec:choice-selection-results} examines how
using these different reference weights affects configuration
selection.

%%%%%%%%%%%%%%%%%%%%%%%%%%%%%%%%%%%%%%%%%%%%%%%%%%%%%%%%%%
\subsubsection{Configuration selection and consequences of robustness (Q4)}
\label{sec:choice-selection-results}

We now address (Q4) by comparing the configurations selected under mean-weight, fitted-weight, and robust optimization. We report decisions for all six combinations of model and role and interpret them alongside the Brier-score improvement $\Delta B$ in~\eqref{eq:choice-prediction}, order agreement, and dominance compliance in Table~\ref{tab:choice-main}.

% Numerical entries retained from sec_4_codex_18set/10_nominal_robust_decisions.py.
\begin{table}[tb]
 \centering
 \caption{Selected configurations under mean (MW), fitted (FW), and uncertain (RB) weights. Entries identify configurations $x^a$. $\Delta W$ (minutes) and $\Delta C$ (weekly opening hours) report RB minus MW; $P$ and $G$ are defined in~\eqref{eq:choice-robust-comparison}. All reported optima are unique.}
 \label{tab:choice-decisions}
 \small
 \setlength{\tabcolsep}{4pt}
 \resizebox{\textwidth}{!}{%
 \begin{tabular}{llccccrrrr}
 \hline
 Model & Role & $\bar\lambda$ & \shortstack{Mean weight\\$x^{\mathrm{MW}}$} & \shortstack{Fitted weight\\$x^{\mathrm{FW}}$} & \shortstack{Robust\\$x^{\mathrm{RW}}$} & $\Delta W$ & $\Delta C$ & $P$ (\%) & $G$ (\%) \\
 \hline
 9B & Patient & 0.750 & E & E & H & +1.483 & -10 & 10.34 & 28.57 \\
 9B & Physician & 0.719 & E & A & H & +1.483 & -10 & 7.70 & 28.57 \\
 9B & Operations & 0.896 & B & A & G & +2.060 & -26 & 31.21 & 22.61 \\
 \hline
 27B & Patient & 0.830 & B & E & E & +1.071 & -18 & 3.63 & 8.61 \\
 27B & Physician & 0.767 & E & E & E & +0.000 & +0 & 0.00 & 0.00 \\
 27B & Operations & 0.662 & E & G & G & +0.989 & -8 & 0.20 & 4.00 \\
 \hline
 \end{tabular}%
 }
\end{table}

Table~\ref{tab:choice-decisions} shows that (MW) and (RB) select different configurations in five of the six combinations. For 27B Patient, (MW) selects $x^{\mathrm B}$, with $(f_w,f_{\mathrm{cost}})=(87.815,74)$, whereas (RB) selects $x^{\mathrm E}$, with $(f_w,f_{\mathrm{cost}})=(88.886,56)$. The robust selection reduces opening hours by 18 per week and increases the average of the group mean waiting times by 1.071 minutes. This change incurs a 3.63\% loss at the mean weight and reduces the worst-case objective by 8.61\%. For 27B Operations, the selection changes from $x^{\mathrm E}$ to $x^{\mathrm G}$, with $(f_w,f_{\mathrm{cost}})=(89.875,48)$. Opening hours decrease by eight and average group waiting time increases by 0.989 minutes, with a 0.20\% loss at the mean weight and a 4.00\% worst-case reduction. Both rules select $x^{\mathrm E}$ for 27B Physician.

For all three 27B roles, (FW) and (RB) select the same configuration: $x^{\mathrm E}$ for Patient and Physician, and $x^{\mathrm G}$ for Operations. Thus, the effect of replacing a single weight by an interval depends on which single-weight rule is used as the reference. For Operations, configurations $x^{\mathrm E}$ and $x^{\mathrm G}$ exchange optimality at $\lambda=0.6583$, which lies between the fitted weight 0.6446 and the bootstrap mean 0.6620. Both weights are close to this switching value, explaining why the (MW) and (FW) selections differ while the loss at the mean weight is only 0.20\%. A small change in the reference weight can therefore change the selected configuration even when its objective consequence is small.

The 9B decisions illustrate the consequences of broad weight uncertainty. For Patient and Physician, $\mathcal U_\lambda=[0,1]$, so~\eqref{eq:choice-robust-endpoints} reduces to $R(x)=\max\{w(x),c(x)\}$. The robust rule selects $x^{\mathrm H}$ for both roles, whereas the mean-weight rule selects $x^{\mathrm E}$. Under this common interval, the robust problem depends only on the normalized performance profiles and retains no information distinguishing the two roles' fitted weights. For Operations, the interval is narrower and the selection changes from $x^{\mathrm B}$ under (MW) to $x^{\mathrm G}$ under (RB). These selections are solutions of the stated optimization problems, but robustness to the chosen weight interval does not resolve the weak predictive and behavioral support for the 9B preference models.

The analysis connects preference inference to configuration selection while retaining the limits of the choice evidence. Outcome differences improve holdout prediction beyond presentation order for two of the six combinations. Accounting for weight uncertainty changes the configuration selected by (MW) in five of the six cases, with higher average group waiting time and fewer MIU opening hours (\textbf{Q4 takeaway}). The associated changes are quantified in waiting minutes, opening hours, and weighted objective values. These are conditional decision consequences for the eleven reported configurations; they do not establish alignment with human stakeholder preferences or clinical effectiveness. Extending the decision set beyond these configurations would require evaluating additional operating schedules and routing decisions through the ED simulation while retaining the chosen preference representation and normalization.

%%%%%%%%%%%%%%%%%%%%%%%%%%%%%%%%%%%%%%%%%%%%%%%%%%%%%%%%%%
\section{Conclusions}\label{sec:future_work}
We proposed LLM-PADA to assess LLM-generated preferences as inputs to multi-objective optimization. The framework considers explicit preferences elicited as criterion importance scores and implicit preferences inferred from pairwise choices between alternatives. It evaluates their reliability, constructs uncertainty sets for the resulting weights, and compares decisions obtained through mean-weight and robust optimization.

In outpatient scheduling, deliberation weakened contextual responses and sometimes reversed role differences. In ED planning, outcome differences improved holdout prediction beyond presentation order for two of six model–role combinations, but choices remained sensitive to order. These results show that agreement and predictive performance alone are insufficient to establish preference reliability.
Compared with mean-weight optimization, the robust models selected different decisions in ten of sixteen outpatient comparisons and five of six ED comparisons. The comparison quantifies the loss at the mean weights incurred to reduce the worst-case weighted objective, for fixed operational data and objective normalizations. This protection applies to the specified weight sets; it does not validate importance scores as objective weights or establish that inferred trade-offs represent stakeholder preferences.

Future work should first compare agents and human stakeholders under the same elicitation tasks, assessing agreement in both the stated trade-offs and the decisions they imply. Replication across model families would test the scope of the findings. Further decision models could explicitly account for uncertainty in simulated outcomes alongside uncertainty in preference weights.

\section*{Acknowledgements}
T. Giovannelli received partial support from the U.S. National
Science Foundation (NSF) [Grant DMR-2555725].

% \bibliography{ref}   % <-- matches ref.bib exactly

% ================================================================
\appendix
% ================================================================

\section{Additional details for explicit preferences}
\label{app:explicit-results}

\subsection{Contextual sensitivity and role differentiation}

Table~\ref{tab:appendix-focal-shifts} reports the paired score differences shown in Figure~\ref{fig:explicit-shift-preservation}. Table~\ref{tab:appendix-role-persistence} gives the complete role-differentiation scores defined in~\eqref{eq:role-persistence}.

\begin{table}[p]
\centering
\caption{Mean paired differences in importance scores between HW and B for waiting and between HO and B for overtime. Pre denotes the initial scores; the remaining columns report results after round 3 of each protocol. Differences are in points on the $1$--$9$ scale, and brackets give 95\% percentile bootstrap intervals over the ten base clinic descriptions.}
\label{tab:appendix-focal-shifts}
\resizebox{\textwidth}{!}{%
\begin{tabular}{lllcccc}
\hline
Condition & Criterion & Agent & Pre & Numerical & Argumentative & Mechanical \\
\hline
HW & Waiting & Intake & $0.560\ [0.225,0.845]$ & $-0.100\ [-0.640,0.440]$ & $0.060\ [-0.320,0.560]$ & $0.580\ [0.384,0.757]$ \\
HW & Waiting & Provider & $0.695\ [0.455,0.910]$ & $0.300\ [-0.280,0.860]$ & $0.040\ [-0.420,0.480]$ & $0.617\ [0.418,0.785]$ \\
HW & Waiting & Scheduler & $0.506\ [0.206,0.736]$ & $0.680\ [0.100,1.320]$ & $0.480\ [0.000,0.960]$ & $0.565\ [0.346,0.738]$ \\
HO & Overtime & Intake & $2.872\ [2.603,3.128]$ & $0.760\ [0.260,1.320]$ & $0.880\ [0.060,1.760]$ & $2.612\ [2.458,2.807]$ \\
HO & Overtime & Provider & $1.691\ [1.167,2.321]$ & $0.980\ [0.140,1.860]$ & $1.600\ [0.660,2.620]$ & $2.287\ [2.009,2.646]$ \\
HO & Overtime & Scheduler & $2.977\ [2.727,3.218]$ & $1.200\ [0.520,1.920]$ & $1.340\ [0.520,2.220]$ & $2.640\ [2.446,2.870]$ \\
\hline
\end{tabular}}
\end{table}

\begin{table}[p]
\centering
\caption{Role-differentiation scores by weight source and condition. Pre denotes scores before deliberation; the other sources use round-3 scores. Each column compares the named agent's normalized weight on its role-associated criterion with the mean weight assigned by the other two agents. Brackets give 95\% percentile bootstrap intervals over the ten base clinic descriptions.}
\label{tab:appendix-role-persistence}
\resizebox{\textwidth}{!}{%
\begin{tabular}{llccc}
\hline
Weight source & Condition & Intake--waiting & Provider--overtime & Scheduler--idle \\
\hline
Pre & B  & $0.087\ [0.071,0.104]$ & $0.094\ [0.079,0.106]$ & $0.078\ [0.064,0.093]$ \\
Pre & CF & $0.142\ [0.129,0.156]$ & $0.127\ [0.118,0.136]$ & $0.087\ [0.075,0.099]$ \\
Pre & HW & $0.073\ [0.044,0.101]$ & $0.091\ [0.071,0.110]$ & $0.094\ [0.068,0.120]$ \\
Pre & HO & $0.053\ [0.035,0.071]$ & $0.019\ [0.010,0.028]$ & $0.026\ [0.012,0.042]$ \\
Numerical & B  & $0.029\ [-0.002,0.060]$ & $-0.039\ [-0.061,-0.013]$ & $0.058\ [0.034,0.080]$ \\
Numerical & CF & $-0.030\ [-0.049,-0.011]$ & $-0.015\ [-0.029,-0.001]$ & $-0.031\ [-0.043,-0.018]$ \\
Numerical & HW & $-0.016\ [-0.038,0.008]$ & $-0.026\ [-0.047,-0.003]$ & $0.042\ [0.023,0.060]$ \\
Numerical & HO & $0.016\ [-0.008,0.041]$ & $-0.031\ [-0.053,-0.011]$ & $0.022\ [0.009,0.035]$ \\
Argumentative & B  & $0.001\ [-0.025,0.023]$ & $-0.019\ [-0.044,0.001]$ & $0.030\ [0.009,0.052]$ \\
Argumentative & CF & $-0.024\ [-0.035,-0.013]$ & $-0.010\ [-0.039,0.018]$ & $-0.022\ [-0.034,-0.009]$ \\
Argumentative & HW & $-0.021\ [-0.043,-0.00003]$ & $0.022\ [-0.011,0.052]$ & $-0.002\ [-0.021,0.018]$ \\
Argumentative & HO & $0.006\ [-0.017,0.032]$ & $0.003\ [-0.014,0.023]$ & $0.004\ [-0.024,0.034]$ \\
Mechanical & B  & $0.023\ [0.019,0.028]$ & $0.026\ [0.022,0.030]$ & $0.022\ [0.018,0.026]$ \\
Mechanical & CF & $0.039\ [0.036,0.043]$ & $0.035\ [0.032,0.037]$ & $0.024\ [0.020,0.027]$ \\
Mechanical & HW & $0.020\ [0.012,0.027]$ & $0.025\ [0.019,0.030]$ & $0.026\ [0.019,0.033]$ \\
Mechanical & HO & $0.014\ [0.010,0.019]$ & $0.005\ [0.003,0.008]$ & $0.007\ [0.003,0.011]$ \\
\hline
\end{tabular}}
\end{table}

\subsection{Preference uncertainty}
\label{app:explicit-uncertainty}

Table~\ref{tab:appendix-uncertainty} reports the componentwise
bounds used to construct the uncertainty sets
$\mathcal U_{S,\pi}$ in~\eqref{eq:explicit-uncertainty-set}.
For each condition and weight source, the lower and upper bounds
are the empirical 5th and 95th percentiles of the normalized
criterion weights. These bounds summarize variation across agents
and base clinic descriptions and, for the numerical and argumentative
protocols, repeated deliberation trajectories. Together with
nonnegativity and the requirement that the three weights sum to one,
they define the weight vectors considered by the robust scheduling
model.

\begin{table}[p]
\centering
\caption{Empirical 5th--95th percentile bounds used with the simplex constraint to define the weight uncertainty sets. For each condition, the pre-deliberation and mechanical samples contain 30 normalized weight vectors, while the numerical and argumentative samples contain 150.}
\label{tab:appendix-uncertainty}
\resizebox{\textwidth}{!}{%
\begin{tabular}{llccc}
\hline
Weight source & Condition & Waiting bounds & Idle bounds & Overtime bounds \\
\hline
Pre & B  & $[0.253,0.400]$ & $[0.236,0.390]$ & $[0.312,0.479]$ \\
Pre & CF & $[0.255,0.471]$ & $[0.172,0.318]$ & $[0.340,0.496]$ \\
Pre & HW & $[0.301,0.456]$ & $[0.204,0.365]$ & $[0.269,0.435]$ \\
Pre & HO & $[0.231,0.327]$ & $[0.176,0.280]$ & $[0.447,0.534]$ \\
Numerical & B  & $[0.272,0.571]$ & $[0.118,0.342]$ & $[0.206,0.469]$ \\
Numerical & CF & $[0.341,0.518]$ & $[0.118,0.267]$ & $[0.300,0.444]$ \\
Numerical & HW & $[0.333,0.615]$ & $[0.125,0.333]$ & $[0.200,0.421]$ \\
Numerical & HO & $[0.294,0.533]$ & $[0.118,0.299]$ & $[0.310,0.471]$ \\
Argumentative & B  & $[0.286,0.615]$ & $[0.074,0.333]$ & $[0.154,0.500]$ \\
Argumentative & CF & $[0.333,0.538]$ & $[0.071,0.267]$ & $[0.308,0.471]$ \\
Argumentative & HW & $[0.355,0.636]$ & $[0.071,0.313]$ & $[0.154,0.471]$ \\
Argumentative & HO & $[0.316,0.554]$ & $[0.071,0.267]$ & $[0.310,0.500]$ \\
Mechanical & B  & $[0.283,0.346]$ & $[0.259,0.353]$ & $[0.318,0.429]$ \\
Mechanical & CF & $[0.325,0.395]$ & $[0.213,0.275]$ & $[0.373,0.427]$ \\
Mechanical & HW & $[0.330,0.419]$ & $[0.252,0.348]$ & $[0.271,0.380]$ \\
Mechanical & HO & $[0.246,0.301]$ & $[0.175,0.268]$ & $[0.455,0.541]$ \\
\hline
\end{tabular}}
\end{table}

\subsection{Complete mean-weight and robust scheduling comparison}
\label{app:explicit-schedules}

Table~\ref{tab:appendix-schedules} compares (MW) and (RB) on the fixed scheduling instance described in Subsection~\ref{sec:explicit-experiment-setup}. The condition identifies the elicitation context from which the weights were obtained; the patient requests, service durations, and scheduling constraints are identical across all rows.

\begin{table}[p]
\centering
\caption{Outcomes selected by (MW) and (RB) on the same scheduling instance, using weights from each source and elicitation condition. Each outcome profile is reported as (total patient waiting, provider idle time, provider overtime) in minutes. The final columns give the percentage changes in the mean-weight and worst-case objectives when moving from the (MW) solution to the (RB) solution, using the corresponding objective value at the (MW) solution as the reference.}
\label{tab:appendix-schedules}
\resizebox{\textwidth}{!}{%
\begin{tabular}{llcccc}
\hline
Weight source & Condition & (MW) & (RB) & Mean-objective change & Worst-case change \\
\hline
Pre & B  & $(86,0,0)$ & $(50,15,6)$ & $+8.12\%$ & $-5.95\%$ \\
Pre & CF & $(72,8,0)$ & $(50,15,6)$ & $+0.99\%$ & $-9.86\%$ \\
Pre & HW & $(72,8,0)$ & $(50,15,6)$ & $+0.01\%$ & $-7.93\%$ \\
Pre & HO & $(86,0,0)$ & $(72,8,0)$ & $+0.16\%$ & $-3.45\%$ \\
Numerical & B  & $(23,25,16)$ & $(23,25,16)$ & $0.00\%$ & $0.00\%$ \\
Numerical & CF & $(50,15,6)$ & $(23,25,16)$ & $+0.57\%$ & $-3.53\%$ \\
Numerical & HW & $(23,25,16)$ & $(23,25,16)$ & $0.00\%$ & $0.00\%$ \\
Numerical & HO & $(50,15,6)$ & $(23,25,16)$ & $+2.20\%$ & $-1.03\%$ \\
Argumentative & B  & $(23,25,16)$ & $(23,25,16)$ & $0.00\%$ & $0.00\%$ \\
Argumentative & CF & $(23,25,16)$ & $(23,25,16)$ & $0.00\%$ & $0.00\%$ \\
Argumentative & HW & $(23,25,16)$ & $(23,25,16)$ & $0.00\%$ & $0.00\%$ \\
Argumentative & HO & $(50,15,6)$ & $(23,25,16)$ & $+3.18\%$ & $-4.37\%$ \\
Mechanical & B  & $(86,0,0)$ & $(86,0,0)$ & $0.00\%$ & $0.00\%$ \\
Mechanical & CF & $(72,8,0)$ & $(50,15,6)$ & $+0.86\%$ & $-3.03\%$ \\
Mechanical & HW & $(72,8,0)$ & $(50,15,6)$ & $+0.12\%$ & $-4.61\%$ \\
Mechanical & HO & $(86,0,0)$ & $(72,8,0)$ & $+0.20\%$ & $-1.18\%$ \\
\hline
\end{tabular}}
\end{table}

\subsection{Prompt templates and validation}
\label{app:explicit-prompts}

The templates below reproduce the wording assembled by the released scripts; angle-bracketed items are substituted programmatically, and long lines are wrapped for typesetting. The same role descriptions are used in elicitation and both LLM deliberation protocols:

\begin{quote}\small
\texttt{Intake:} You represent the Intake/Patient perspective. You focus on patient experience and timely access while recognizing that clinic capacity and provider workload also matter.\par
\texttt{Provider:} You represent the Provider/Physician perspective. You focus on a sustainable clinical workload, preserving clinical quality, and limiting work that extends beyond the scheduled end of the session.\par
\texttt{Scheduler:} You represent the Scheduler/Operations perspective. You focus on feasible and efficient use of available clinic capacity, reducing unused appointment gaps, and maintaining an orderly schedule.
\end{quote}

The elicitation prompt has the following form. Each call requests all three scores together, with one of the preceding role descriptions inserted into \texttt{<ROLE>}:

\begin{small}
\begin{verbatim}
Assign importance scores for an outpatient scheduling decision.

Role: <ROLE>

<OPERATIONAL STATE>

Criteria:
- waiting: patient waiting time before service starts
- idle: provider idle time between patients
- overtime: provider work beyond the scheduled session end

Score all three criteria jointly on a 1-9 scale, where 1 is very low
importance and 9 is extremely high importance. Preserve your stakeholder
perspective, but respond to the operational state when one is shown.
Use the scale comparatively and do not make every criterion equally high.
Return only a JSON object with exactly these integer fields: waiting,
idle, overtime.

\end{verbatim}
\end{small}

For B, HW, and HO, \texttt{<OPERATIONAL STATE>} contains the nine saved fields in this order: \texttt{current\_time}, \texttt{session\_end}, \texttt{patients\_seen}, \texttt{patients\_remaining}, \texttt{open\_slots}, \texttt{mean\_waiting\_time}, \texttt{provider\_idle\_time}, \texttt{projected\_overtime}, and \texttt{cancellations}. For CF, the entire block is exactly \texttt{Operational state: not provided (context-free condition).} No stored B values are shown.

The prompts request criterion importance; they do not present the seven-request scheduling instance, candidate appointment patterns, or the reference scales of 120, 60, and 60 minutes later used in the objective. Using the normalized importance scores as weights for these objectives is an assumption of the scheduling experiment. The elicitation wording for idle time refers to intervals between patients; the scheduling implementation additionally counts any initial gap before the first service, as specified above.

Both deliberation prompts begin with the same text:

\begin{small}\begin{verbatim}
Update criterion importance scores in a structured outpatient scheduling
deliberation.

Round: <ROUND> of 3
Role: <ROLE>

<OPERATIONAL STATE>

Current score vectors for all stakeholders (1-9 scale):
<JSON OBJECT: Intake, Provider, Scheduler; each with waiting, idle,
overtime>

Revise waiting, idle, and overtime jointly after considering the other
vectors. Remain faithful to your stakeholder role and the shown
operational state. Do not move merely to create consensus. Scores must
be integers from 1 to 9.
\end{verbatim}
\end{small}

Numerical deliberation appends exactly:

\begin{small}\begin{verbatim}
Return only JSON with exactly: waiting, idle, overtime.
Do not provide a reason.
\end{verbatim}
\end{small}

Argumentative deliberation instead appends a prior-argument block followed by the instruction below. At round 1 the block reads \texttt{No prior-round arguments exist. referenced\_argument must be NONE.} At later rounds, it is a JSON array of the three preceding argument records, each containing \texttt{argument\_id}, \texttt{stakeholder}, the three applied scores, \texttt{reason}, \texttt{objection}, and \texttt{parse\_status}.

\begin{small}\begin{verbatim}
In addition to the same three scores, give a one-sentence reason of
at most 30 words, identify one other stakeholder whose position you
address, give a one-sentence objection of at most 30 words, and list
only closed state-variable identifiers actually used. In CF use
['none']. If you rely on an immediately preceding argument, cite its
exact ID; otherwise use NONE.
Allowed state_variables_used: <STATE-VARIABLE IDENTIFIERS>.
Return only the specified JSON object.
\end{verbatim}
\end{small}

The argumentative response requires \texttt{waiting}, \texttt{idle}, \texttt{overtime}, \texttt{state\_variables\_used}, \texttt{reason}, \texttt{target\_agent}, \texttt{referenced\_argument}, and \texttt{objection}. Ollama receives a strict JSON schema in addition to the text prompt. The scripts independently validate the responses and reject self-targeting and invalid argument references. Each deliberation request permits up to two retries. If no valid response is obtained, the failure is recorded and the agent's previous-round score vector is carried forward. Every raw attempt is retained. The complete instantiated prompts and outputs are available in the JSONL audit files.

\end{document}